\documentclass[smallextended,nospthms]{svjour3}       
\smartqed  
\usepackage{amssymb,amsmath,amsfonts,amsthm,url}
\usepackage[usenames,dvipsnames]{color}
 \usepackage{lipsum}

\newtheorem{Theorem}{Theorem}[section]
\newtheorem{Lemma}[Theorem]{Lemma}
\newtheorem{prop}[Theorem]{Proposition}
\newtheorem{Remark}[Theorem]{Remark}
\newtheorem{Example}[Theorem]{Example} 
\newtheorem{cor}[Theorem]{Corollary}
\newtheorem{Assumption}{Assumption}[section]
\newtheorem{Definition}[Theorem]{Definition}

\theoremstyle{definition}
\numberwithin{equation}{section}

\def\oF{\overline F}
\def\la{\leftarrow}

\def\conv{\stackrel{v}{\to}}

\def\a{\alpha}
\def\b{\beta}

\def\E{\mathbb{E}}

\def\M{\mathbb{M}}
\def\P{\mathbb{P}}
\def\R{\mathbb{R}}

\def\bX{\boldsymbol X}
\def\bY{\boldsymbol Y}

\def\ba{\boldsymbol a}

\def\bu{\boldsymbol u}
\def\bx{\boldsymbol x}
\def\by{\boldsymbol y}

\def\bzero{\boldsymbol 0}

\def\binfty{\boldsymbol \infty}

\def\RV{\mathcal{RV}}
\def\MRV{\mathcal{MRV}}
\def\VaR{\text{VaR}}

\definecolor{darkred}{RGB}{139,0,0}
\definecolor{darkgreen}{RGB}{0,139,0}

\allowdisplaybreaks 
\numberwithin{equation}{section}
\begin{document}


\title{Risk concentration under second order regular variation
\thanks{Bikramjit Das gratefully acknowledges partial support from  MOE2017-
T2-2-161. Partial support from RARE-318984 (an FP7 Marie Curie IRSES Fellowship) is kindly acknowledged by both authors. 
}}


\author{Bikramjit Das         \and
       Marie Kratz }


\institute{B. Das \at
            Singapore University of Technology and Design\\
            8 Somapah Road, Singapore 487372\\
              \email{bikram@sutd.edu.sg}          
           \and
          M. Kratz \at
             ESSEC Business School, CREAR\\ 
             avenue Bernard Hirsch BP 50105\\ 
             Cergy-Pontoise 95021 Cedex, France\\
             \email{kratz@essec.edu}     
}

\date{Received: date / Accepted: date}

\titlerunning{Risk concentration}
\authorrunning{Das and Kratz}
\maketitle

\begin{abstract}
Measures of risk concentration and their asymptotic behavior for portfolios with heavy-tailed risk factors is of interest in risk management. Second order regular variation is a structural assumption often imposed on such risk factors to study their convergence rates.
In this paper, we provide the asymptotic rate of convergence of the \emph{measure of risk concentration} for a portfolio of heavy-tailed risk factors, when the portfolio admits the so-called \emph{second order regular variation} property. 
Moreover, we explore the relationship between multivariate second order regular variation for a vector (e.g.,  risk factors) and the second order regular variation property for the sum of its components (e.g., the portfolio of risk factors). 
Results are illustrated with a variety of examples. 

\end{abstract}


\section{Introduction}\label{sec:intro}

An important issue in risk management is assessing the effects of adding an investment to a portfolio of risk factors (time series of returns) and understanding how this aggregate risk relates to the individual risk factors. Broadly studied under the labels of  \emph{risk concentration} or \emph{risk diversification}, the past couple of decades have seen tremendous developments in the understanding of  this topic. Our interest is in a portfolio of risk factors that are heavy-tailed, where adequate care is necessary to study the aggregation of the risk factors; see \cite{dacorogna:elbahtouri:kratz:2018,embrechts:mcneil:straumann:2002,ibragimov_etal:2011,kratz:2014,puccetti:ruschendorf:2013}  for detailed discussions on diversification, especially under heavy-tailed returns. 

We concentrate on the particular risk measure {\it{value-at-risk}}. Recall that for a random variable (risk factor) $X$ with distribution function $F$, the value-at-risk  at level $0<\beta<1$ is defined as
\begin{align*}
\VaR_{\b}(X) : = \inf \{y\in \R: \P(X\le y) \ge \beta\} = F^{\la}(\b).
\end{align*}
Consider a portfolio of risk factors $\bX=(X_{1},\ldots, X_{d})$. Unless otherwise specified, we assume for this paper that $X_{1},\ldots, X_{d}$ are  identically distributed (homogeneous) non-negative random variables. The behavior of the sum $$S_{d} = X_{1}+ \ldots + X_{d}$$ and its value-at-risk $\VaR_{\beta}(S_{d})$ have been studied under various assumptions, either on the marginal distribution $F$ (where $X_{1}\sim F$) or on the dependence structure of $\bX$.
If $X_{1}, \ldots, X_{d}$ are independent and identically distributed (iid) with a regularly varying tail distribution having tail parameter $\alpha>0$, then it is well-known that $\VaR(S_{d})$ is asymptotically sub-additive or super-additive according as $\alpha>1$ or $\alpha<1$ (see \cite{degen_etal:2010,embrechts:lambrigger:wuthrich:2009}); an accurate estimation for high threshold has been proposed in \cite{kratz:2014}. 
Further refinements of such first order approximations have been studied under the notion of second order regular variation for the marginal tail $\overline{F}$; see \cite{degen_etal:2010,dehaan:ferreira:2006,kortschak:2012,mao:hu:2013,resnickbook:2008}. 
Furthermore, there are studies on the asymptotic behavior of the tail of $S_{d}$ and $\VaR_{\b}(S_{d})$ under specific copula assumptions on the dependence structure of $\bX$ (see \cite{alink_etal:2004,barbe_etal:2006,kortschak:2012,sun:li:2010}), or by providing risk bounds under assumptions on marginal densities (see \cite{peng:wang:yang:2013,puccetti:ruschendorf:2013}).

We observe that an assumption of second order regular variation is often useful in finding rate of convergence of various risk measures in the heavy-tailed regime. Hence an investigation of univariate and multivariate second order regular variation and their interplay is important. Thus the broad goal of this paper is twofold. First, we  assess the limit behavior of the measure of risk concentration for value-at-risk of the sum of risk factors under the assumption that the sum is second order regular varying. Secondly, we identify conditions under which multivariate second order regular variation for a vector would imply second order regular variation property for the sum of its components; thus the rate of convergence of risk concentration measures could be obtained in such cases.

The paper is structured as follows. In Section \ref{sec:regvar}, we discuss the various notions of regular variation both first order and second order as well as univariate and multivariate, and posit an alternative definition of univariate second order regular variation in Proposition~\ref{prop:alt2rv}. 
In Section \ref{sec:divindex}, by assuming second order regular variation for the sum of the components of a random vector $\bX\in \R_+^d$, we find the convergence rate of the diversification benefit $D_{\beta}(\bX)$ for the $\VaR$ risk measure as $\beta\to 1$.
We explore the effect of multivariate second order regular variation ($2\MRV$) on aggregation in Section~\ref{sec:2mrvto2rv}. Here we provide sufficient conditions on a $2\MRV$ vector for the sum of its components to satisfy a second order regular variation ($2\RV$) property.
We illustrate our results with examples from popularly used marginal distributions and copula dependence structures in Section \ref{sec:ex}. We also show through an example that $\bX\in2\MRV$ is not a necessary condition for the sum of its components to be $2\RV$.  Conclusions are drawn in Section \ref{sec:concl}. In the Appendix, the diversification property relating the marginal risks to the aggregate risk is extended to tail equivalent risks. 

\subsection{Notation} \label{subsec:notation}

A summary of some notation and concepts used in this paper is provided here.  We use bold
  letters to denote vectors, with capital letters for random vectors
  and small letters for non-random vectors, e.g.,
  $\by=(y_1,\ldots,y_d)\in \R^d$. We also define
  $\bzero=(0,\ldots,0)$ and $\binfty=(\infty,\ldots,\infty)$. Vector operations are always understood component-wise,
  e.g., for vectors $\bx, \by \in \R^d$, $\bx\le \by$ means $x_i\le y_i$
  for  $1\le i\le d$.  For a constant $k\in \R$ and a set $A\subset \R^{d}$, we denote by $kA:= \{k\bx: \bx\in A\}$.
Some additional notation follows.
Detailed discussions are in the references provided.

\begin{tabular}{p{1.3cm}p{12.7cm}ll}
\\[2mm]
  $\mathbb{E}$ & $[0,\infty]^{d}\setminus \{\bzero\}$\\[2mm]
  $\mathcal{B} (\mathbb{E})$ & The Borel $\sigma$-field of the subspace $\mathbb{E}$.\\[2mm]
  $\mathbb{M}_+(\mathbb{E})$ & The class of Radon measures on Borel subsets of $\mathbb{E}$.\\[2mm]
      $\conv$ & vague convergence of measures, often on  $\mathbb{M}_+(\mathbb{E})$; see \cite{resnickbook:2007}.\\[2mm]
            $\aleph$ & The set $\{\bx\in \E: ||\bx||=1\}$, where $||\cdot||$ denotes the Euclidean norm in $\R^{d}$.\\[2mm]
%
          $f^{\leftarrow}$ & The left-continuous inverse of a monotone function $f$.\\[2mm]
  $\RV_{\rho}$ & The class of regularly varying functions (at $\infty$) with index $\rho\in \R$,  that
         is, functions $f:\mathbb{R}_+\to \mathbb{R}_+$
 satisfying $\displaystyle \lim_{t\to\infty}f(tx)/f(t)=x^\rho,$ for $x>0$; see  \cite{bingham:goldie:teugels:1989,dehaan:ferreira:2006,resnickbook:2008}.\\[2mm]
$f\sim g$ & Functions $f:\mathbb{R}_+\to \mathbb{R}_+$ and $g:\mathbb{R}_+\to \mathbb{R}_+$ are \emph{asymptotically equivalent} or \emph{tail equivalent} if $\displaystyle{\lim_{t\to\infty} f(t)/g(t)=1}$ and we  write $f\sim g$.
\end{tabular}


\section{Preliminaries on Regular Variation}\label{sec:regvar}

Regular variation often forms the basis for studying heavy-tailed distributions. In this section, we recall definitions and properties of first and second order regular variation in both univariate and multivariate cases \cite{bingham:goldie:teugels:1989,dehaan:1970,dehaan:ferreira:2006,resnick:2002,resnickbook:2008}. 
Furthermore, we introduce an alternative definition of univariate second order regular variation in  Proposition~\ref{prop:alt2rv}, which exhibits a form of the limit measure that is often obtained in practice.

\subsection{Regular variation in one dimension} \label{subsec:regvar-dim1}

\begin{Definition}[Regular variation, \cite{bingham:goldie:teugels:1989}]\label{def:rv}
A random variable $X$ with distribution function $F$ has regularly varying (right) tail with index $\alpha > 0$ if $\oF = 1-F \in \RV_{-\alpha}$.
Alternatively, we  say that $X$ has a regularly varying tail if there exists a function $b:\R_{+}\to\R_{+}$ with $b(t) \uparrow \infty$ as $t\to \infty$ such that
\begin{align*}  
& \lim_{t\to\infty} t\,\P[X>b(t) x] = x^{-\a} .  
\end{align*}
We write $\oF\in \RV_{-\alpha}(b)$ or, by abuse of notation, $X \in \RV_{-\alpha}(b)$.

\end{Definition}

A consequence of the definition is that $b(\cdot)\in \RV_{1/\a}$ and a natural choice is $b(t)=(1/\overline{F})^{\la}(t)$. For example, Burr, Fr\'echet, Pareto, stable distributions with shape parameter $\alpha\in(0,2)$, have $\RV_{-\alpha}$ tail distributions (see  e.g. \cite[p. 35]{embrechts:kluppelberg:mikosch:1997} for further details).  Furthermore, some distribution functions $F$, with $\overline{F} \in \RV_{-\alpha}$,  have a second order property that is not captured by the  scaling in the definition of regular variation. The Pareto-Lomax distribution, analyzed in Example \ref{ex:parlomax}, is one such distribution and the following definition describes this property. 

\begin{Definition}[Second order regular variation;  \cite{dehaan:resnick:1993}, \cite{resnick:2002}]\label{def:2rv}
A random variable $X$ with distribution function $F$ such that $\oF\in \RV_{-\alpha}(b)$ with 
$\alpha > 0$, possesses second order regular variation with parameter $\rho\le 0$, 
if there exists a function $\displaystyle A(t)\underset{t\to\infty}{\to} 0$ that is ultimately of constant sign,  $|A(\cdot)|\in \RV_{\rho}$ with $\rho\le 0$ {and $c\neq 0$} such that
\begin{equation}\label{eq:def2RV}
\frac{t\oF(b(t)x)-x^{-\alpha}}{A(b(t))} \underset{t\to\infty}{\longrightarrow} cx^{-\alpha} \frac{x^{\rho}-1}{\rho} =: H(x), \quad x>0.
\end{equation}
The right hand side of \eqref{eq:def2RV} is interpreted as $H(x)=cx^{-\alpha}\log(x)$ when $\rho=0$. We write $\oF \in 2\RV_{-\alpha,\rho}(b,A,H)$ or, by abuse of notation, $X\in 2\RV_{-\alpha,\rho}(b,A,H)$. Some or all of the arguments within the brackets may be omitted for simplicity.
\end{Definition}

\begin{Remark} \label{rk:equivResDeHaan} \normalfont
An equivalent representation of \eqref{eq:def2RV} suggests (see \cite{dehaan:ferreira:2006}): $\oF \in 2\RV_{-\alpha,\rho}(A^{*},H^{*})$ if there exists an ultimately positive or negative function $A^*$ with $\displaystyle A^*(t)\underset{t\to\infty}{\to} 0$ such that 
\begin{align}\label{eq:def2RValt}
 \lim_{t\to\infty}\frac{\frac{\oF(tx)}{\oF(t)}-x^{-\alpha}}{A^{*}(t)} = H^{*}(x) := cx^{-\alpha} \frac{x^{\rho}-1}{\rho} , \quad x>0,
 \end{align}
for some constant $c\neq 0$ and parameters $\alpha>0, \rho\le 0$. 
The parameters $\a,\rho$ of course remain the same in both definitions \eqref{eq:def2RV} and  \eqref{eq:def2RValt}. By replacing $t$ with  $b(t)=\left(1/{\oF}\right)^{\la}\!\!\!(t)$, the functions $A^{*}$ and $A$ are asymptotically equivalent, and $H^{*}$, $H$ coincide.
\end{Remark}

\begin{Remark}\label{rem:Hnondeg}
As a consequence of Theorem 1 in \cite{dehaan:stadtmueller:1996},
 it is known that if \eqref{eq:def2RValt} holds with $H^{*}(x)$ which is not a multiple of $x^{{-\alpha}}$, then $H^{*}$ necessarily satisfies the given representation. In our proofs we will often use \eqref{eq:def2RValt} as a definition of $2\RV$.\end{Remark}

\begin{Remark}\label{rk:2rvunifconv}
The convergences \eqref{eq:def2RV} and \eqref{eq:def2RValt} in the definition of second order regular variation  hold locally uniformly on compact intervals of $(0,\infty)$; cf.  \cite{geluk:dehaan:1987,dehaan:stadtmueller:1996}. This fact is particularly useful in proving some of the convergence results that we obtain.
\end{Remark}

\begin{Example}\label{ex:parlomax}
Consider the Pareto-Lomax distribution function given by $\overline F(x)=(1+x)^{-\alpha}, x>0$ and $\alpha>0$. Choosing $b(t)=(1/\overline{F})^{\la}(t)=t^{1/\a}-1$ and $A(t)=(1+t)^{-1}$, we obtain
$$
\lim_{t\to \infty}\frac{t\oF(b(t)x) - x^{-\a}}{A(b(t))} 
= -\a x^{-\a}(x^{-1}-1)=:H(x), \quad x>0.
$$
Hence $\oF \in 2\RV_{-\a,-1}(b,A,H)$ (with $c=\a$ in \eqref{eq:def2RV}).
\end{Example}

From \cite{dehaan:ferreira:2006}, we know that the function of $t$ in \eqref{eq:def2RValt} has an asymptotic limit $H^{*}$ necessarily of the form given. An important consequence of Definition \ref{def:2rv} is that there may be other tail equivalent choices of $A$ and $b$ that would provide a non-zero limit in \eqref{eq:def2RV}. The form of the limit measure in this case is modified a bit.  Nevertheless we may always get back the original definition by appropriate choices of functions $b$ and $A$. 
This idea helps in finding the appropriate constants for the limit quantities obtained in Theorem \ref{thm:2mrvtoSrv}. The following proposition formalizes the idea.

\begin{prop}\label{prop:alt2rv}
Let $X$ be a random variable with  distribution function $F$ such that $\oF\in \RV_{-\alpha}(b)$, $\alpha > 0$, and assume there exists a function $\displaystyle A(t)\underset{t\to\infty}{\to} 0$ that is ultimately of constant sign,  $|A(\cdot)|\in \RV_{\rho}$ with $\rho\le 0$, and $c\neq 0, c^* \in \R$ with $c\neq c^*\rho$, such that
\begin{align}\label{eq:def2RVinprop}
\frac{t\oF(b(t)x)-x^{-\alpha}}{A(b(t))} &\underset{t\to\infty}{\longrightarrow} cx^{-\alpha} \frac{x^{\rho}-1}{\rho} + c^*x^{-\alpha}=:H_{(c,c^*)}(x), \quad x>0.
\end{align}
Then there exist functions $\tilde b \in \RV_{1/\a}$ and  $\tilde A(t)\underset{t\to\infty}{\to} 0$ that is ultimately of constant sign and $|\tilde A| \in \RV_{\rho}$, such that 
\begin{equation}\label{eq:def2RVinprop-bis}
X\in 2\RV_{-\a,\rho}(\tilde b, \tilde A, \tilde H)\quad\text{where} \;\tilde H(x)=cx^{-\a} \frac{x^{\rho}-1}{\rho}, \;  \quad x>0. 
\end{equation}
Note that when $\rho=0$, the RHS of \eqref{eq:def2RVinprop} is $cx^{{-\alpha}}\log(x)+c^*x^{-\alpha}$.
\end{prop}

\begin{proof}
Let $\rho<0$. Since $\overline{F}\in \RV(b)$, we have $b(t) \uparrow \infty$ as $t\to\infty$, (w.l.o.g. it can be chosen strictly increasing). Hence its inverse,  $b^{\la}(s): = \inf\{t\in\R: b(t)\ge s\}$ is also a strictly increasing function, so that $b(b^{\la}(s))= s = b^{\la}(b(s)).$
Moreover $\displaystyle \lim_{t\to \infty} A(t)\to 0$. Hence, as $s\to\infty$, 
\[\tilde{b}(s):= \left(\frac{b^{\la}(\cdot)}{1+ c^*  A(\cdot)}\right)^{\la}(s)\, \to \infty.\]
 Now w.l.o.g. assume that $\tilde{b}:(0,\infty)\to(0,\infty)$ is also strictly increasing (otherwise we may choose an appropriate $A^*\sim A$).
Also note that $\tilde{b}(t)\sim b(t)$ as $t\to\infty$.   Moreover, let $\tilde{A}\equiv A$. Then we have, for any fixed $x>0$,
 \begin{align*}
\frac{t\oF(\tilde{b}(t)x)-x^{-\alpha}}{\tilde{A}(\tilde{b}(t))} &= \frac{\tilde{b}^{\la}(s)\overline{F}(sx)-x^{-\alpha}}{A(s)} \quad\quad (\text{where }\; s:=\tilde{b}(t))\\
&  \!\!\!\!  \!\!\!\!\!\!\!\! \sim\frac{{b}^{\la}(s)\overline{F}(sx)-(1+c^* A(s))x^{-\alpha}}{(1+c^* A(s))A(s)}
 =   \frac{{b}^{\la}(s)\overline{F}(sx)-x^{-\alpha}}{A(s)} \cdot \frac{1}{1+c^* A(s)}  -   \frac{1}{1+c^* A(s)} c^* x^{-\alpha}\\
& \!\!\!\!  \!\!\!\!\!\!\!\! =  \frac{t'\overline{F}(b(t')x)-x^{-\alpha}}{A(b(t'))} \cdot \frac{1}{1+c^* A(b(t'))}  -  \frac{1}{1+c^* A(s)} c^* x^{-\alpha} \quad (\text{where }\; t'={b}^{\la}(s)={b}^{\la}(\tilde{b}(t)).\\
\end{align*}
Letting $t \to \infty$ (which implies $t'\to \infty$), we deduce from this last equality that
 \begin{align*}
\lim_{t\to\infty} \frac{t\oF(\tilde{b}(t)x)-x^{-\alpha}}{\tilde{A}(\tilde{b}(t))} 
& = cx^{-\alpha} \frac{x^{\rho}-1}{\rho} + c^*x^{-\alpha} -c^*x^{-\alpha} = cx^{-\alpha} \frac{x^{\rho}-1}{\rho}.
\end{align*}
We can show the case where $\rho=0$ in a similar manner.
\end{proof}
\vspace{1ex}
\begin{Remark}\label{rem:other2rv}
Henceforth we may use \eqref{eq:def2RVinprop} to mean that $X\in2\RV_{-\a,\rho}(b,A,H_{(c,c^*)})$ as an alternative to the limit identities given in \eqref{eq:def2RV} or \eqref{eq:def2RValt}. Note that, for $\rho<0$, we can write
\begin{align}\label{eq:def2RVHcc}
H_{(c,c^*)}(x) = \tilde{c}x^{-\alpha}\frac{Mx^\rho-1}{\rho}
\end{align}
where $\tilde{c} =c-c^*\rho\neq 0$ and $M=c/\tilde{c}$; this is also a form observed for $2\RV$ sometimes. \\
Since $c\neq c^{*}\rho$, the RHS of \eqref{eq:def2RVinprop} is not a multiple of $x^{-\alpha}$, which implies that $X\in 2\RV$. 
\end{Remark}

\subsection{Regular variation in dimension $d>1$} \label{subsec:regvar-dim2}

Multivariate regular variation facilitates the study of jointly heavy-tailed random variables and is a natural extension of Definition \ref{def:rv}.
The following definitions explain multivariate regular variation as well as second order regular variation for joint tail distributions of random variables. The notion of vague convergence of measures is used for convergence of measures on the non-negative Euclidean orthant $\R_{+}^{d}$ and its subsets; see \cite{resnickbook:2007} for further details.

\begin{Definition}[Multivariate regular variation, \cite{resnickbook:2007}]\label{def:mrv}
 Suppose $\bX=(X_1,\ldots,X_d)$ is a random vector in  $[0,\infty)^d$. Then $\bX$ is {multivariate regularly varying}  if there exist $b(t) \uparrow \infty$ with $b(\cdot) \in \RV_{1/\alpha}$, $\a>0$, and a Radon measure $\nu\neq 0$ such that as $t\to\infty$
\begin{align*}
t\,\P\left(\frac{\bX}{b(t)} \in \;\cdot\; \right) \, \conv  \, \nu(\cdot) \quad \text{on } \;\; \M_{+}(\E),
\end{align*}
where $\conv $ denotes vague convergence of measures on the space $\M_{+}(\E)$. We write  $\bX\in \MRV_{-\alpha}(b,\nu)$.
\end{Definition}
The measure $\nu(\cdot)$ has a scaling property for relatively compact $A\subset \E$, given by
\begin{align*}
\nu(kA) = k^{-\alpha}\nu(A), \quad k>0\,.
\end{align*}
We restrict attention to univariate second order regular variation for the first part of the paper. Eventually we exhibit connections between second order regular variation and multivariate regular one that we define below.

\begin{Definition}[Second order multivariate regular variation, \cite{resnick:2002}]\label{def:2mrv} 
Suppose $\bX\in \MRV_{-\alpha}(b,\nu)$ and there exists $\displaystyle A(t)\underset{t\to\infty}{\to} 0$ that is ultimately of constant sign with $|A(\cdot)| \in \RV_{\rho},\, \rho \le 0$,  such that
\begin{align}\label{eq:2mrv}
\frac{t\,\P\left(\frac{\bX}{b(t)} \in [\bzero,\bx]^c \right) - \nu([\bzero,\bx]^c)}{A(b(t))} \underset{t\to\infty}{\to} H(\bx)
\end{align}
locally uniformly in $\bx\in (0,\infty]^d \setminus \{\binfty\}$, where $H$ is a function that is non-zero and  finite. Then $\bX$ is {\it{second order regularly varying}} with parameters $\alpha > 0$ and $\rho\le 0$. We write {$\bX\in 2\MRV_{-\alpha,\rho}(b,A,\nu,H)$; some or all of the parameters may be omitted according to the context.}
\end{Definition}
\begin{Remark} A regularly varying distribution with tail parameter $\alpha=0$ is called super-heavy-tailed and we avoid this case for the purposes of this paper.
\end{Remark}
\begin{Remark}
 Replacing $d=1$ in Definitions~\ref{def:mrv}  and \ref{def:2mrv} gives back the univariate versions, i.e., Definitions ~\ref{def:rv} and ~\ref{def:2rv} respectively. In order to use \eqref{eq:2mrv} in terms of vague convergence of signed measures, we need to impose further conditions on the distribution $F$ of $\bX$ as aptly noted in \cite[Section 4]{resnick:2002}.
\end{Remark}


\section{Measuring Risk Concentration} \label{sec:divindex}

\subsection{Risk measures and the diversification index}

In risk management, evaluating  risk concentration (or, equivalently,  diversification benefit) properly is key for both insurance and investments. 
Let $\mathcal{L}$ be the collection of random variables defined on a given probability space $(\Omega,\mathcal{F},\P)$. For a risk measure $\eta:\mathcal{L}\to\R \cup\{-\infty,\infty\}$, and  $d$ risks $X_i\in \mathcal{L}, i=1,\cdots,d$,
 the associated  {\it measure of risk concentration} or {\it{diversification index}} (see \cite{buergi_etal:2008,tasche:2008})  for $\bX=(X_{1},\ldots,X_{d})$, is given by
\begin{align*}
 D_\eta(\bX) & = \frac{\eta\left(\sum\limits_{i=1}^d X_i\right)}{\sum\limits_{i=1}^d \eta(X_i)}.
  \end{align*}

We understand that the  risk measure is sub-additive, additive  or super-additive according as $D_{\eta}$ less than, equal to or more than 1 respectively. In this work, we concentrate on the  popular risk measure \emph{value-at-risk} (VaR) as the choice for $\eta$ and obtain asymptotic results for $D_{\VaR_{\beta}}$. Moreover, we assume that all univariate marginals of $\bX$ are identically distributed, if not otherwise specified.

Properties of the diversification index $D_{\VaR_{\beta}}$ and its asymptotic limits under different assumptions on the marginal distributions and dependence structure, can be found in the literature (see e.g. \cite{buergi_etal:2008,dacorogna:elbahtouri:kratz:2018,degen_etal:2010,embrechts:kluppelberg:mikosch:1997}). Henceforth, we denote  $D_{\VaR_{\b}}$ as $D_{\b}$ emphasizing its dependence on $\b$. The following result is illustrative and we study its extension under a more general set-up.

\begin{Lemma}[see \cite{embrechts:lambrigger:wuthrich:2009}, Example 3.1]\label{lem:divIndex-iid}
Suppose $\bX=(X_1,\ldots, X_d) \in \R_+^d$ has iid components with \linebreak $X_{1} \in \RV_{-\a}$, $\alpha>0$. Let $S_d :=\sum_{i=1}^d X_i$. Then
\begin{equation}\label{eq:limDivI-iidRV}
\lim_{\b\uparrow 1} D_{\b}(\bX) = \lim_{\b\uparrow 1}\frac{\VaR_\b (S_d)}{d \; \VaR_\b (X_1)}  = d^{1/\a-1}.
\end{equation}
\end{Lemma}

The precise rate of convergence  for the limit in \eqref{eq:limDivI-iidRV} as $\beta\to 1$, can be obtained by using an additional assumption of second order regular variation (marginally) on $\oF$; see \cite{degen_etal:2010,mao:hu:2013,omey:willekens:1986}. Some studies relax the condition of independence of marginal variables and obtain limits as in \eqref{eq:limDivI-iidRV} as well as rates of convergence; for instance \cite{hua:joe:2011b} works under a scale-mixture dependence with  second order regularly varying marginal distributions, \cite{kortschak:2012} works under an assumption of  asymptotic independence, and \cite{tong:wu:xu:2012} assumes an Archimedean copula as the dependence structure. Note that the result in \cite{tong:wu:xu:2012} for Archimedean copulas is not entirely correct, but we find comparable results for Clayton copula in Example \ref{ex:parclay-lom}.

\vspace{5pt}
{Prior to stating our main result, we state an auxiliary result on the effect of assuming that a random variable $X\in 2\RV$, on its  value-at-risk $\VaR(X)$. Results in a similar spirit to this can also be found in the literature, for example in \cite[Section 2.3]{dehaan:ferreira:2006}.}

\begin{Lemma}\label{lem:applicVervaat}
For any positive random variable $X \in 2\RV_{-\a,\rho}(b,A,H)$ where $H(x)=\frac{c}{\rho}x^{-\alpha}(x^{\rho}-1)$ with $c\neq 0$, we have
\begin{equation*}\label{eq:VaR-for-2RV}
  \lim_{t\to\infty}\frac{1}{A(b(t))}\left(\frac{\VaR_{1-x/t}(X)}{b(t)}- x^{-1/\alpha}\right)  = H^*(x) = \begin{cases}  
\frac{c}{\a\rho}x^{-1/\a}(x^{-\rho/\a}-1) & \text{if}\; \rho<0,\\[0.5em]
- \frac{c}{\a^{2}}\log(x) &\text{if}\; \rho=0.
\end{cases}
\end{equation*}
The above convergence also holds locally uniformly on $(0,\infty)$.
\end{Lemma}

\begin{proof}[Proof of Lemma \ref{lem:applicVervaat}]
The proof is given for $\rho<0$; the case for $\rho=0$ can be done in a similar way. We use \emph{Vervaat's Lemma} (see \cite{vervaat:71}, \cite[Lemma A.0.2]{dehaan:ferreira:2006}), which is stated here for notational convenience. \\
{ \sc{Vervaat's Lemma:}}
\emph{Suppose that $y$ is a continuous function on $[0,\infty)$, $\{z_t(.)\}_{t\ge0}$ is a family of non-negative, non-increasing functions, and $g$ is a function having a negative continuous derivative. Let $\delta(t) \to 0$ with $\delta(t)>0$ eventually, such that
$$
\lim\limits_{t\to\infty} \frac{z_t(x) - g(x)}{\delta(t)} = y(x),
$$
locally uniformly on $(0,\infty)$. Then, locally uniformly on $(g(\infty), g(0))$,
$$
\lim\limits_{t\to\infty} \frac{z_t^{\leftarrow}(x) - g^{\leftarrow}(x)}{\delta(t)} = -\left(g^{\leftarrow}\right)'(x)\;y(g^{\leftarrow}(x))\,.
$$
}
Since $X \in 2\RV_{-\a,\rho}(b,A,H)$, $A(\cdot)$ is of constant sign ultimately. W.l.o.g. assume $A(t)>0$ eventually (otherwise we choose $-A(t)$ as the $A(t)$).  Now denoting  $z_t(x) = t\,\P[X>x b(t)]=t\overline F_X(xb(t))$, $g(x) = x^{-\alpha}$, $\delta(t) = A(b(t))$ and $y(x) = H(x)$, and applying Vervaat's Lemma to \eqref{def:2rv} (which holds locally uniformly on $(0,\infty)$), we get, for $x>0$:
\begin{align*}
 \lim_{t\to\infty}\frac{1}{A(b(t))}\left(\frac{\VaR_{1-x/t}(X)}{b(t)}- x^{-1/\alpha}\right)  &  =   \lim_{t\to \infty} \frac{\frac{1}{b(t)} \oF_{X}^{\leftarrow} (x/t) - x^{-1/\alpha}}{A(b(t))} = \lim_{t\to\infty} \frac{z_t^{\leftarrow}(x)-g^{\leftarrow}(x)}{\delta(t)} \\
 & = -\left(g^{\leftarrow}\right)'(x)\;y(g^{\leftarrow}(x)) = \frac{c}{\a\rho}x^{-1/\a}(x^{-\rho/\a}-1)= H^*(x),
\end{align*}
and the convergence holds locally uniformly.
\end{proof}

\begin{Example}\label{ex:vervaat}
The following example provides a simple application of Lemma \ref{lem:applicVervaat}.  \\
Suppose $X\sim F$, with 
$\displaystyle \overline{F}(x)= 1-F(x)= (x^{-\a}+x^{-2\a})/2, \quad x>1.$
Hence, for $0<p<1$, $\overline{F}^{\la}(p) = 2^{1/\a}(\sqrt{1+8p}-1)^{-1/\a}$. With $b(t) = \overline{F}^{\la}(1/t)$, we have, for $x>0$,
\begin{align*}
 t\overline{F}(b(t)x) & = \frac t2 \left[\frac 12\left(\sqrt{1+\frac 8t}-1\right)x^{-\a} +\frac 14\left(\sqrt{1+\frac 8t}-1\right)^{2}\!\! x^{-2\a}\right]=x^{-\a} \left(1+\frac 2t(x^{-\a}-1) + o(t^{-1})\right)\\
&\underset{t\to\infty}{\to}x^{-\a}.
 \end{align*}
Moreover, taking $A(t)=t^{-\a}$, we obtain $A(b(t)) = 2\left[1-2/t + o(1/t)\right]/t$, from which we deduce that
\[ 
\lim_{t\to \infty} \frac{t\overline{F}(b(t)x) - x^{-\a}}{A(b(t))} = x^{-\a}(x^{-\a}-1) =: H(x).
\]
Hence,  $X \in 2\RV_{-\alpha,-\alpha}(b,A,H)$ with $c=-\a$ as defined in \eqref{eq:def2RV}. 
Applying Lemma \ref{lem:applicVervaat},  we have, for $x>0$,
\[  \lim_{t\to\infty}\frac{1}{A(b(t))}\left(\frac{\VaR_{1-x/t}(X)}{b(t)}- x^{-1/\alpha}\right)   =  \frac1{\a} \,x^{-1/\a}(x-1), \]
which may also be directly verified. 
\end{Example}

\subsection{The diversification index of $\VaR_{\beta}$ under second order regular variation for the sum}

We are concerned with the limit and the convergence rate of the diversification benefit \linebreak $D_{\beta}(\bX) = D_{{\VaR_{\beta}}}(\bX)$ for $\beta$ close to 1. We restrict to random vectors $\bX\in \R_+^d$ with the sum of the components denoted $S_d = X_1+ \ldots +X_d$. 
Our main result Theorem~\ref{thm:riskconc} concerns the asymptotic rate of convergence for the diversification index $D_{\beta}(\bX)$ under a 2$\RV$ assumption on the sum and uses Lemma~\ref{lem:applicVervaat} to obtain rates of convergences for value-at-risk. Theorem~\ref{thm:riskconc}(1) is a generalization of Lemma \ref{lem:divIndex-iid} relaxing the iid assumption and has been observed previously (see \cite{barbe_etal:2006,degen_etal:2010}); we state it for the sake of completeness. Theorem~\ref{thm:riskconc}(2) provides the precise rate of convergence under  assumptions both on $S_{d}$ and $X_{1}$. In Section \ref{sec:2mrvto2rv}, we provide conditions and examples under which such a result is applicable.

\begin{Theorem}\label{thm:riskconc} 
Let $\bX\in \R_+^d$ be such that $\bX\in \MRV_{-\alpha}(b,\nu)$ with identical marginal distributions.
\begin{enumerate}
\item[1.] Then we have
$$
 \lim_{\b\uparrow 1}  D_{\b}(\bX) = \frac 1d\left(\frac{\nu(\Gamma_d)}{\nu(L_1)}\right)^{1/\alpha}  =:K_d
$$
where
\begin{equation}\label{def:Gamma_d}
\Gamma_d= \{\bx\in \R_+^d: x_1+\ldots+x_d>1\} \quad \text{and} \quad  L_{1}=\{\bx\in \R_+^d: x_1>1\}.
\end{equation}
\item[2.] Assume that $S_d\in 2\RV_{-\alpha,\rho}(b_d,A_d,H_d)$ where $\displaystyle{H_d(x):= \frac{c_d}{\rho}x^{-\alpha}(x^{\rho}-1)}$, $x>0$, and $\rho<0$. 
\begin{enumerate}
\item[a.] Suppose $X_1\in 2\RV_{-\alpha,\rho_1}(b_{1},A_1,H_1)$ with $\displaystyle{H_1(x):= \frac{c_1}{\rho_{1}}x^{-\alpha}(x^{\rho_{1}}-1)}$, $x>0$, and $\rho_{1}<0$. 
\begin{enumerate}
\item[i.] If $\displaystyle{\lim_{t\to\infty}A_1(b_{1}(t))/A_d(b_{d}(t)) =\kappa\in \R}$, then, for any $x>0$, we have
\begin{align}\label{lim:D1minxtk}
  \lim_{t \to \infty} \frac{D_{1- x/t}(\bX) - K_{d} }{A_{d}(b_{d}(t))}\, =\,
\frac{(c_{d}-\kappa c_{1})K_{d}}{\alpha\rho} (x^{{-\rho/\a}}-1).
\end{align}
\item[ii.] If $\displaystyle{\lim_{t\to\infty}A_d(b_{d}(t))/A_1(b_{1}(t)) =0}$, then, for any $x>0$, we have
\begin{align}\label{lim:D1minxtnok}
  \lim_{t \to \infty} \frac{D_{1- x/t}(\bX) - K_{d} }{A_{1}(b_{1}(t))}\, =\,
-\frac{c_{1}K_{d}}{\alpha\rho_{1}} (x^{{-\rho_1/\a}}-1).
\end{align}
\end{enumerate}
\item[b.]  Suppose $X_1$ does not possess  $2\RV$ and $\displaystyle{\frac{\VaR_{1-x/t}(X_{1})}{b_{}(t)} - x^{-1/\alpha} (\nu(L_1))^{1/\alpha} = o(A_{d}(b_{d}(t)))}$  as $t\to\infty$, then, for any $x>0$,
\begin{align}\label{lim:D1minxtno2rv}
  \lim_{t \to \infty} \frac{D_{1- x/t}(\bX) - K_{d} }{A_{d}(b_{d}(t))}\, =\,
\frac{c_{d}K_{d}}{\alpha\rho} (x^{{-\rho/\a}}-1).
\end{align}
\end{enumerate}
\end{enumerate}

\end{Theorem}
 \newpage

\begin{Remark}
\normalfont
~\vspace{1pt}
\begin{enumerate}
\item[(i)] Theorem  \ref{thm:riskconc}(2) provides sufficient conditions to obtain the rate of convergence of $D_{\beta}$. Merely assuming $S_{d} \in 2\RV$ does not provide enough information. 
 If the marginal variable $X_{1}\in 2\RV$ as observed in part (2.a), we may obtain different rates of convergence depending on the behavior of the  auxiliary functions $A_{1}(t)$ vis-a-vis $A_{d}(t)$. On the other hand, if $X_{1}$ does not possess $2\RV$, we can still obtain a rate of convergence using part (2.b).
\item[(ii)]  If $c_{d}=\kappa c_{1}$, then part (2.a.i) of Theorem \ref{thm:riskconc} is clearly not informative enough. Otherwise, we always obtain non-zero limits in \eqref{lim:D1minxtk}-\eqref{lim:D1minxtno2rv} since $c_{1}, c_{d}, K_{d},\alpha,\rho$ are all non-zero.
\item[(iii)] The rate of convergence can also be found when either one or both of $\rho$ and $\rho_{1}$ are zero, using Lemma \ref{lem:applicVervaat} in the appropriate case. We omit these cases to avoid notational confusion. 
\item[(iv)] Statement 2.a. of Theorem~\ref{thm:riskconc} can be extended to the case of tail equivalent risks \linebreak
$\displaystyle X\in 2\RV_{-\alpha,\rho_{ X}}(b_{ X},A_{ X},H_{ X})$ and $\displaystyle Y\in 2\RV_{-\alpha,\rho_{ Y}}(b_{ Y},A_{ Y},H_{ Y})$ such that $\displaystyle 0<\lim_{t\to\infty} \frac{b_X(t)}{b_Y(t)}<\infty$; see Theorem~\ref{theo:appendix}  in the Appendix.
\end{enumerate}
\end{Remark}

\begin{proof}  [Proof of Theorem \ref{thm:riskconc}]
~\vspace{1pt}
\begin{itemize}
\item[1.]Since $\bX\in \MRV_{-\alpha}(b,\nu)$, for $x>0$, we have
\begin{align}
 \lim_{t\to\infty}  t\,\P(X_1>xb(t)) & =  \lim_{t\to\infty}  t\,\P\left(\frac{\bX}{b(t)} \in x L_{1} \right) = x^{-\alpha}\nu(L_1), \label{lam1} \\
 \text{and}\quad \lim_{t\to\infty}  t\,\P(S_d>xb(t)) &=  \lim_{t\to\infty}  t\,\P\left(\frac{\bX}{b(t)} \in x \Gamma_{d}\right) = x^{-\alpha}\nu(\Gamma_d). \label{gamd}
  \end{align}
 Inverting \eqref{lam1} and \eqref{gamd} (see \cite[Proposition A.0.1]{resnickbook:2008}), we obtain, for any $x>0$,
\begin{align}\label{lamgam}
\lim_{t\to\infty}\frac{\VaR_{1-x/t}(X_1)}{b(t)} = x^{-1/\alpha}(\nu(L_1))^{1/\alpha} \quad \text{and} \quad \lim_{t\to\infty}\frac{\VaR_{1-x/t}(S_d)}{b(t)} = x^{-1/\alpha}(\nu(\Gamma_d))^{1/\alpha}.
\end{align}
Hence we have 
\begin{align*}
\lim_{\b\uparrow 1}  D_{\b}(\bX) = \lim_{t\to\infty} D_{1-1/t} (\bX) 
 & =  \lim_{t\to\infty}   \frac 1d \cdot \frac{\VaR_{1-1/t}(S_d)}{b(t)} \cdot \frac{b(t)}{\VaR_{1-1/t}(X_1)}
   = \frac 1d\left(\frac{\nu(\Gamma_d)}{\nu(L_1)}\right)^{1/\alpha} \!\!\!\!\!\!=K_{d}.
\end{align*}
\item[2.] 
Now $S_{d} \in \RV_{{-\alpha}}$ is evident from \eqref{gamd}. By assuming $S_d\in 2\RV_{-\alpha,\rho}(b_d,A_d,H_d)$, we also get
   $ \displaystyle{\lim_{t\to\infty}  t\;\P(S_d>b_d(t)) = 1}$, which, compared with \eqref{gamd}, implies that $b_d(t) \underset{t\to\infty}{\sim} (\nu(\Gamma_d))^{1/\alpha} b(t)$. 
W.l.og. we assume $b_d(t) = (\nu(\Gamma_d))^{1/\alpha} b(t)$.  Since $S_d\in 2\RV_{-\alpha,\rho}(b_d,A_d,H_d)$, applying  Lemma~\ref{lem:applicVervaat} gives, for $x>0$,
\begin{align}\label{lim:2mrvSd}
&\lim_{t\to\infty}\frac1{A_{d}(b_{d}(t))} \left(\frac{\VaR_{1-x/t}(S_{d})}{b_{d}(t)}- x^{-1/\alpha}\right) = \frac{c_{d}}{\alpha \rho} x^{-1/\alpha}(x^{-\rho/\alpha}-1)=: H_d^*(x).
\end{align}
To assess the second order property, observe that, for any $x>0$,
$$
 {D_{1- x/t}(\bX) - K_{d} }  = \frac{\VaR_{1-x/t}(S_{d})}{d\VaR_{1-x/t}(X_{1})}- K_{d}  =: {\bf{I}}(x,t) - {\bf{II}}(x,t),
$$
where $\displaystyle {\bf{I}}(x,t):= \frac{b_d(t)}{d\VaR_{1-x/t}(X_{1})} \left[\frac{\VaR_{1-x/t}(S_{d})}{b_{d}(t)} - x^{-1/\alpha}\right]$\\
 and $\displaystyle {\bf{II}}(x,t):=\frac{K_{d}\,b(t)}{\VaR_{1-x/t}(X_{1})} \left[\frac{\VaR_{1-x/t}(X_{1})}{b(t)} - x^{-1/\alpha} (\nu(L_1))^{1/\alpha}\right]$.\\
Now, using  \eqref{lamgam} and \eqref{lim:2mrvSd}, we have for any $x>0$,
\begin{align}\label{boldI}
 \lim_{t\to \infty} \frac{{\bf{I}}(x,t)}{A_{d}(b_{d}(t))} = K_d\cdot x^{1/\a} \cdot H_{d}^{*}(x) = K_{d}\frac{c_{d}}{\alpha \rho}(x^{-\rho/\alpha}-1).
\end{align}
For analyzing ${\bf{II}}(x,t)$, using \eqref{lamgam}, we know that for any $x>0$,
\[ \lim_{t\to\infty}\frac{K_{d}\,b(t)}{\VaR_{1-x/t}(X_{1})} = {K_{d}\, x^{1/\alpha}(\nu(L_{1}))^{{-1/\alpha}}} >0.  \]
{\emph{Case (a):}} We have assumed $X_1\in 2\RV_{-\alpha,\rho_1}(b_{1},A_1,H_1)$. Using Lemma \ref{lem:applicVervaat}, we obtain for $x>0$,
\begin{align}\label{x1h1}
\lim_{t\to\infty}\frac{1}{A_{1}(b_{1}(t))} \left[\frac{\VaR_{1-x/t}(X_{1})}{b_{1}(t)} - x^{-1/\alpha} \right] =  \frac{c_1}{\alpha \rho_{1}} x^{-1/\alpha}(x^{-\rho_{1}/\alpha}-1)=: H_1^*(x).
\end{align}
W.l.o.g., we assume $b_{1}(t) =(\nu(L_{1}))^{1/\alpha} b(t)$ (cf. \eqref{lamgam}). 

\begin{quote}
{\emph{Sub-case (a.i.):}}  Since we further assume $\displaystyle{\lim_{t\to\infty}A_1(b_{1}(t))/A_d(b_{d}(t)) =\kappa\in \R}$, we can write
\begin{align}
 \lim_{t\to \infty} \frac{{\bf{II}}(x,t)}{A_{d}(b_{d}(t))}  & =  \lim_{t\to \infty} \frac{1}{A_{d}(b_{d}(t))}\times\frac{K_{d}b(t)}{\VaR_{1-x/t}(X_{1})} \times\left[\frac{\VaR_{1-x/t}(X_{1})}{b(t)} - x^{-1/\alpha} (\nu(L_1))^{1/\alpha}\right] \nonumber\\
 & =  \lim_{t\to \infty}   \frac{A_{1}(b_{1}(t))}{A_{d}(b_{d}(t))}\times\frac{K_{d}b_{1}(t)}{\VaR_{1-x/t}(X_{1})}\times  \frac{1}{A_{1}(b_{1}(t))} \left[\frac{\VaR_{1-x/t}(X_{1})}{b_{1}(t)} - x^{-1/\alpha} \right] \nonumber\\
  &= \kappa \times K_{d}x^{1/\alpha} \times H_{1}^{*}(x), 
\end{align}
which is 0 only if $\kappa=0$. Note that, if $\kappa\neq 0$, we can conclude that $\rho=\rho_{1}$.  Hence, 
\begin{align}
 \lim_{t\to \infty} \frac{{\bf{II}}(x,t)}{A_{d}(b_{d}(t))}  & =   \frac{ \kappa c_{1}K_{d}}{\alpha \rho} (x^{-\rho/\alpha}-1). \label{boldII2}
\end{align}
Therefore, using \eqref{boldI} and \eqref{boldII2}, we have, for any $x>0$,
\begin{align*}
 \lim_{t\to \infty} \frac{D_{1-x/t}-K_{d}}{A_{d}(b_{d}(t))} =  \lim_{t\to \infty}\left[ \frac{{\bf{I}}(x,t)}{A_{d}(b_{d}(t))} - \frac{{\bf{II}}(x,t)}{A_{d}(b_{d}(t))} \right]= \frac{ (c_{d}-\kappa c_{1})K_{d}}{\alpha \rho}(x^{-\rho/\alpha}-1).
\end{align*}
\indent {\emph{Sub-case (a.ii.):}}  In contrast to the previous part, we assume  $\displaystyle{\lim_{t\to\infty}A_d(b_{d}(t))/A_1(b_{1}(t)) =0}$. Hence, using \eqref{boldI}, we have
 \begin{align}\label{boldInew}
 \lim_{t\to \infty} \frac{{\bf{I}}(x,t)}{A_{1}(b_{1}(t))} =  \lim_{t\to \infty} \frac{A_d(b_{d}(t))}{A_1(b_{1}(t))}\times  \frac{{\bf{I}}(x,t)}{A_{d}(b_{d}(t))} =0.
\end{align}
On the other hand, from \eqref{lamgam} and \eqref{x1h1}, we obtain
\begin{align}
 \lim_{t\to \infty} \frac{{\bf{II}}(x,t)}{A_{1}(b_{1}(t))}  & =  \lim_{t\to \infty} \frac{K_{d}b(t)}{\VaR_{1-x/t}(X_{1})}  \times \frac{1}{A_{1}(b_{1}(t))} \left[\frac{\VaR_{1-x/t}(X_{1})}{b(t)} - x^{-1/\alpha} (\nu(L_1))^{1/\alpha}\right] \nonumber\\
  &= K_{d}x^{1/\alpha} \times H_{1}^{*}(x) 
   =  \frac{ c_{1}K_{d}}{\alpha \rho_{1}} (x^{-\rho_{1}/\alpha}-1). \label{boldII2new}
\end{align}
We deduce,  using \eqref{boldInew} and \eqref{boldII2new}, that, for any $x>0$,
\begin{align*}
 \lim_{t\to \infty} \frac{D_{1-x/t}-K_{d}}{A_{1}(b_{1}(t))} =  \lim_{t\to \infty}\left[ \frac{{\bf{I}}(x,t)}{A_{1}(b_{1}(t))} - \frac{{\bf{II}}(x,t)}{A_{1}(b_{1}(t))} \right]= -\frac{c_{1}K_{d}}{\alpha \rho_{1}}(x^{-\rho_{1}/\alpha}-1).
\end{align*}
\end{quote}
{\emph{Case (b):}} If $X_{1}$ does not possess $2\RV$ and
$\displaystyle{\frac{\VaR_{1-x/t}(X_{1})}{b(t)} - x^{-1/\alpha} (\nu(L_1))^{1/\alpha} = o(A_{d}(b_{d}(t)))},$
then clearly
\begin{align}\label{boldII1}
 \lim_{t\to \infty} \frac{{\bf{II}}(x,t)}{A_{d}(b_{d}(t))} = 0.
\end{align}
Therefore, combining \eqref{boldI} and \eqref{boldII1} provides, for any $x>0$,
\begin{align*}
 \lim_{t\to \infty} \frac{D_{1-x/t}-K_{d}}{A_{d}(b_{d}(t))} =  \lim_{t\to \infty}\left[ \frac{{\bf{I}}(x,t)}{A_{d}(b_{d}(t))} - \frac{{\bf{II}}(x,t)}{A_{d}(b_{d}(t))} \right]= \frac{c_{d}K_{d}}{\alpha \rho}(x^{-\rho/\alpha}-1).
\end{align*}
\end{itemize}
\vspace{-3ex}
\end{proof}

Now a proportional growth rate of $D_{1-1/t}(X)$ can be deduced immediately from Theorem \ref{thm:riskconc} giving us  the following corollary. 
\begin{cor}\label{cor:dbeta}
Under the conditions of Theorem \ref{thm:riskconc}, we have, for any $x>0, y>0$,
$$
\lim_{t\to\infty} \frac{D_{1-x/t}(\bX)-D_{1-1/t}(\bX)}{D_{1-y/t}(\bX)-D_{1-1/t}(\bX)} = \frac{x^{-\rho^{*}/\a}-1}{y^{-\rho^{*}/\a}-1},
$$
for some $\rho^{*}\le0$. If either $X_1\in 2\RV_{-\alpha,\rho_1}(b_{1},A_1,H_1)$ with  $\displaystyle{\lim_{t\to\infty}A_1(b_{1}(t))/A_d(b_{d}(t)) =\kappa\in \R}$, \\or $X_{1}$ does not possess $2\RV$ but $\displaystyle{\frac{\VaR_{1-x/t}(X_{1})}{b(t)} - x^{-1/\alpha} (\nu(L_1))^{1/\alpha} = o(A_{d}(b_{d}(t)))}$, then $\rho^{*}=\rho$. \\ 
On the other hand, if $X_1\in 2\RV_{-\alpha,\rho_1}(b_{1},A_1,H_1)$ with  $\displaystyle{\lim_{t\to\infty}A_d(b_{d}(t))/A_1(b_{1}(t)) =0}$, then $\rho^{*}=\rho_{1}$.
\end{cor}

Under the assumption that we can statistically estimate $D_{\beta}$ at moderately high values of $\beta$, Corollary~\ref{cor:dbeta} may provide a way to extrapolate values of $D_{\beta}$ to extreme levels of  $\beta$. For instance, suppose our data allows us to compute   estimates of the diversification index for $\VaR$ at $90\%$ and $95\%$, which is given by  $\hat{D}_{0.90}(\bX)$ and $\hat{D}_{0.95}(\bX)$, then for any $p\gg 0.95$ (with $0<p<1$), we may use Corollary~\ref{cor:dbeta} to estimate  $D_{p}(\bX)$ as
\[
\hat{D}_{p}(\bX) = \hat{D}_{0.90}(\bX) + \left(\frac{\left(\frac{1-p}{0.1}\right)^{-\rho^*/\alpha}-1}{\left(0.5\right)^{-\rho^*/\alpha}-1}\right) \left[\hat{D}_{0.95}(\bX) - \hat{D}_{0.90}(\bX)\right].
\]

\begin{Example} \label{subsec:Exappend} 
\normalfont

We illustrate an application of Theorem \ref{thm:riskconc} in the following example. Consider $\bX=(X_{1},X_{2})$ 
with identical marginal distribution $F$ defined by  
$$
\displaystyle  \overline{F}(x):=1-F(x)=\left(1+x\right)^{-\alpha}, \quad\forall x>0.
$$  
Recall, from Example \ref{ex:parlomax}, that $X_1 \in 2\RV_{-\alpha,-1}(b_1,A_1,H_1)$ where $b_1(t)=t^{1/\alpha}-1, A_1(t)=(1+t)^{-1}$ and $H_1(x)=-c_1x^{-\alpha}(x^{-1}-1)$ with $c_1=\alpha$. The dependence structure of $\bX$  is assumed to be a Clayton copula with parameter $\theta>0$,  given by
\[
C_{\theta}\left(u,v\right)=\left(u^{-\theta} + v^{-\theta} - 1 \right)^{-1/\theta}, \quad u,v\in [0,1].
\]
First, we observe that $\bX\in \MRV_{-\alpha}(b,\nu)$, where $b:=b_1$, since
\begin{equation}\label{valueNu1}
t\;\P\left(\frac{\bX}{b(t)}  \in ([0,x_1]\times[0,x_2])^c\right) 
\underset{t\to\infty}{\to}  x_1^{-\a} + x_2^{-\a} - \left(x_1^{\alpha\theta} +x_2^{\alpha\theta} \right)^{-1/\theta} =:  \nu\left( ([0,x_1]\times[0,x_2])^c\right).
\end{equation}
 Let $\theta=1/\alpha$. In this case, $\nu\left( ([0,x_1]\times[0,x_2])^c\right):=x_1^{-\a} + x_2^{-\a} - \left(x_1+x_2\right)^{-\a}$, and,
for $\alpha>1$, using results from \cite[Proposition 2.1]{dacorogna:elbahtouri:kratz:2018}, we find that
\begin{align*} 
 \P(S_2>x) & = \overline F_{2}(x)= \int_{x}^{\infty}\alpha(\alpha+1) \frac{y}{(1+y)^{\a+2}}\,  \mathrm d y 
                          = \frac1{(1+x)^\a}\left(1+ \frac{\a\,x}{1+x} \right),
\end{align*}
where $S_{2}=X_{1}+X_{2}$. Therefore
$$\displaystyle
\lim_{t\to\infty}\frac{\overline F_{2} (tx)}{\overline F_{2}(t)} \;=\; x^{-\a} \qquad \text{and}  \qquad
\lim_{t\to\infty}\frac{\frac{\oF_{2} (tx)}{\overline F_{2}(t)} - x^{-\a}}{\frac{\a(2+\a)}{1+\a}\,t^{-1}} \;= \; - x^{-\a} \left(x^{-1}-1\right).
$$
Hence we have $S_2\in 2\RV_{-\a, -1}(A^*,H^*)$ with  
$$ 
A^*(t):= \,\frac{\a(2+\a)}{1+\a}\,t^{-1}  \quad \text{and}  \quad H^*(x):=- c\,x^{-\a} \left(x^{-1}-1\right)
$$ 
with $c=1$ (see \eqref{eq:def2RValt}). From Remark \ref{rk:equivResDeHaan}, we can also say that
\;$\displaystyle S_2\in 2\RV_{-\a, -1}(b_2,A_2,H_2)$, with 
$$ b_{2}(t) = (1/\overline{F_2})^{\leftarrow}(t)\underset{t\to\infty}{\sim} (\alpha+1)^{1/\alpha}t^{1/\alpha}, \quad  A_2= A^*, \quad H_2= H^*$$ where the constant $c_2=1$ in $H_2$, i.e.,
\begin{equation}\label{eq:ex3.7}
\frac{t\,\P(S_2/b_2(t)>x) - x^{-\alpha}}{\a(\a+2)(\a+1)^{-(1+1/\a)}t^{-1/\a}} \;\underset{t\to\infty}{\to}\;  -x^{-\a}(x^{-1}-1).
\end{equation}
Note that 
$$ \kappa :=\lim_{t\to\infty} \frac{A_1(b_1(t))}{A_2(b_2(t))} =  \frac{(\a+1)^{1+1/\alpha}}{\a(\a+2)}>0.$$
We need to  compute $\displaystyle  K_2:=\lim_{\beta\to1} D_{\beta}$. Using the density function of $\nu$ at $x_1>0,x_2>0$,
$$
\lambda(x_1,x_2)  = \Bigg| \frac{\partial^2}{\partial x_1\partial x_2}  \nu\left( ([0,x_1]\times[0,x_2])^c\right)\Bigg| = \a(\a+1)(x_1+x_2)^{-(\a+2)} \!, 
$$
we can find $\displaystyle \nu(L_1)=\nu\left( ([0,1]\times[0,\infty))^c\right) =1$ and 
\begin{equation}\label{eq:nugamma2}
\nu(\Gamma_2)= \int_{\Gamma_2}  \lambda(x_{1}, x_{2}) \,\mathrm dx_{1}\,\mathrm d x_{2} = \a(\a+1)\int_{\Gamma_2} (x_1+x_2)^{-(\a+2)}\, \mathrm dx_1\,\mathrm dx_2 \, = \, 1+\a.
\end{equation}
Therefore
$$K_2 = \frac 12\left(\frac{\nu(\Gamma_2)}{\nu(L_1)}\right)^{1/\alpha} =\frac 12 (1+\a)^{1/\a}.$$
Hence, applying   Theorem \ref{thm:riskconc}, case (2.a.i.), we obtain, for any $x>0$,
\begin{equation}\label{eq:DVex3.7}
\lim_{t\to\infty} \frac{D_{1-x/t} -  K_2 }{A_2(b_2(t))} = \frac{(1-\kappa \a)}{-\a} (x^{1/\a}-1)  
%
= \frac1{\a}\left(\frac{(\a+1)^{1+1/\a}}{\a+2}-1 \right) (x^{1/\a}-1).
\end{equation}
\end{Example}

\section{Effect of Multivariate Second Order Regular Variation on Aggregation} \label{sec:2mrvto2rv}

In this section, we discuss the relations between multivariate second order regular variation ($2\MRV$) for a vector and the $2\RV$ property for the sum of its components. We look for conditions on a  $2\MRV$ vector $\bX$  to deduce the $2\RV$ property for the sum, illustrating the main result with examples (covering all possible cases considered in the theorem). Then we question if $2\MRV$ is a necessary condition for the sum to be $2\RV$. We show that it is not, providing examples.

\subsection{Main result}\label{subsec:aggr}

\subsubsection{Assumptions}

Here we present the framework for stating the main result. 
Second order regular variation for vector valued random entities has been appropriately discussed in \cite{resnick:2002}.  
 In the following, we provide conditions under which the second order regular variation condition of Definition \ref{def:2mrv} can be represented as vague convergence of measure, depending on whether the limit measure $\nu(\cdot)$ as obtained in Definition \ref{def:mrv} has a density with respect to the Lebesgue measure or not.
Assumption \ref{assumption1} gives the appropriate conditions when $\nu(\cdot)$ has a density with respect to the Lebesgue measure;
note that $\bX$ is then not asymptotically independent. 
When $\nu(\cdot)$ does not have a density, conditions are given in Assumption \ref{assumption2}. As an example, if the tail distribution of $\bX$ exhibits asymptotic independence, then $\nu(\cdot)$ does not have a density.
Let us present the two assumptions on $\bX$, a $d$-dimensional non-negative random vector with distribution function $F$.

\begin{Assumption} \label{assumption1} (see \cite{resnick:2002}, Section 4.2)
\begin{enumerate}
\item  Let $F$ have a density $F'$  and identical one-dimensional marginals $F_1$ such that $\overline{F}_1 \in \RV_{-\alpha}$. Assume that, for $\bx\in \E$,
\begin{align}\label{cond4}
\lim_{t \to\infty} \left|\frac{F^{'}(t\bx)}{t^{-d}\overline{F}_{1}(t)}- \lambda(\bx)\right| = 0 \quad \text{and} \quad \lim_{t \to\infty} \sup\limits_{\ba\in\aleph}\left|\frac{F^{'}(t\ba)}{t^{-d}\overline{F}_{1}(t)}- \lambda(\ba)\right| = 0,
\end{align}
where $\lambda(\cdot) \neq 0$ is bounded on $\aleph$. The limit function $\lambda$ necessarily satisfies $\lambda(t\bx)=t^{{-\alpha-d}} \lambda(\bx)$. 
 \item Assume that the second order condition given in \eqref{eq:def2RV} holds for $\overline{F}_{1}$, so that $\overline{F}_{1}\in 2\RV_{-\alpha,\rho}(A)$,  
with $\rho \le 0$ 
and, for $\bx\in \E$,
\begin{align}\label{cond5}
\lim_{t \to\infty} \left|\frac{\frac{F^{'}(t\bx)}{t^{-d}\overline{F}_{1}(t)} -\lambda(\bx)}{A(t)} -\chi' (\bx)\right| = 0 \quad \text{and} \quad   \lim_{t \to\infty}  \sup\limits_{\ba\in\aleph} \left|\frac{\frac{F^{'}(t\ba)}{t^{-d}\overline{F}_{1}(t)} -\lambda(\ba)}{A(t)} -\chi'(\ba)\right| = 0,
\end{align}
where $\chi'\neq 0$ is integrable on sets bounded away from $\bzero$, and is finite and bounded on $\aleph$.
\end{enumerate}
\end{Assumption}

\begin{Remark}\label{rem:assumption1}
Let $\nu$ be defined as $ \nu([\bzero,\bx]^{c})=\!\!\!\int\limits_{[\bzero,\bx]^{c}} \!\!\!
 \lambda(\bu) \mathrm d\bu$, $\bx>\bzero$.  We define the signed measure by
\begin{equation}\label{signmeas}
\mu_{t}([\bzero,\bx]^{c}) := \frac{t\,\P\left[\frac{\bX}{b(t)} \in [\bzero,\bx]^{c}\right] - \nu ([\bzero,\bx]^{c})}{A(b(t))}
\end{equation}
that has a density given by  \;$\displaystyle \mu'_{t}([\bzero,\bx]^{c}) := \frac{tb(t)^{d}F'(b(t)\bx)- \lambda(\bx)}{A(b(t))}$, $\bx\in [0,\infty)^{d}$.
\end{Remark}

\begin{Assumption} \label{assumption2} (see \cite{resnick:2002}, Section 4.2)
\begin{enumerate}
\item Suppose \eqref{eq:2mrv} holds with $\nu([\bzero,\bx]^{c}) = \kappa\sum_{i=1}^{d} x_{i}^{-\alpha}$, where $\kappa$ is some constant.
\item Assume that the one-dimensional marginals $F_1$ are identical and satisfy the second order condition as in Definition \ref{def:2rv} ({\it i.e.} $\overline F_1\in 2\RV_{-\alpha,\rho}(b,A,{\chi_{1}})$ with $\alpha>0$ and $\rho<0$), such that we also have
\begin{align}\label{cond7}
\mu_{t1}^{\pm} (\cdot):= \left(\frac{t\,\P\left[\frac{X_{1}}{b(t)}\in \cdot \right] -\nu_{\alpha}(\cdot)}{A(b(t))}\right)^{\pm} \conv \chi_{1}^{\pm} (\cdot),\quad \text{as} \; t\to\infty,
\end{align}
on $(0,\infty]$, where $\nu_{\alpha}(x,\infty)=x^{-\alpha}, x>0$ and for $t>0$, $\mu_{t1}^{+}, \mu_{t1}^{-}, \chi^{+}_1, \chi^{-}_1$ are positive Radon measures with $\mu_{t1}=\mu_{t1}^{+}-\mu_{t1}^{-}$ and $\chi_{1}=\chi_{1}^{+}-\chi_{1}^{-}$.
\end{enumerate}
\end{Assumption}

In order to  aggregate multiple risks factors (with the same marginal distribution or at least  equivalent tail order), multivariate regular variation helps in providing justification for sub- or super-additivity; see \cite{degen:embrechts:2011,embrechts:lambrigger:wuthrich:2009}.  We observe here that further structure and intuition can be provided by assuming a second order  regular variation condition. Here we provide a connection between multivariate second order regular variation of $\bX$ and second order regular variation of sum of its components. This eventually helps us in evaluating risk measures for sums of homogeneous random factors with different dependence structures.

\subsubsection{Main result}

Aggregation of risk under multivariate regular variation is relatively straightforward to check. For example, assuming that $\bX\in \MRV_{-\alpha}(b)$ with identical marginal distributions $F_1$,
we can check that, if  $S_d:=\sum_{i=1}^d X_i \sim F_{S_{d}}$ for $d\ge 2$, then $F_{S_{d}} \in \RV_{-\alpha}$ with the same function $b(\cdot)$ as in Definition \ref{def:mrv}. 
 The following proposition extends this implication  to the case where $\bX$ and its components  possess second order regular variation.
 \begin{Theorem}\label{thm:2mrvtoSrv} 
Let $\bX \in \R_{+}^{d}$ with distribution function $F$ such that  $\bX \in 2\MRV_{-\alpha,\rho}(b,A,\nu,H)$. Let $\chi:\mathcal{B}(\E)\to\R$ be such that $\displaystyle \chi([\bzero,\bx]^c):=H(\bx)$ and 
$z\mapsto \chi(z\Gamma_d)$ is neither zero, nor a multiple of $z^{-\alpha}$,  $\Gamma_d$ being defined in \eqref{def:Gamma_d}.
Also assume that $F$  satisfies either Assumption \ref{assumption1} or Assumption \ref{assumption2} (in both assumptions, $\bX$ has identical marginal distributions that are $2\RV$, and to fix notations let  $X_1\in 2\RV_{-\alpha,\rho}{(b,A,H_1)}$ with $H_{1}=\frac{c_1}{\rho}(x^{\rho}-1)$). 
 Then we have the following.
\begin{enumerate}
\item[(1)] The sum $S_d$ satisfies
\begin{equation}\label{eq:bAHstar}
S_d \in 2\RV_{-\alpha, \rho}(b^*,A^*,H^*), \quad \text{where}
\end{equation}
$b^*(t):= (\nu(\Gamma_{d}))^{1/\alpha} b(t),\;\,
A^*(t):=A((\nu(\Gamma_d))^{-1/\alpha}t)$ and 
$H^*(x) := \chi(x(\nu(\Gamma_d))^{1/\alpha} \Gamma_d) = H_{c_a,c_b}(x)$ as defined in  \eqref{eq:def2RVinprop} for some $c_a\neq 0, c_b\in\R$. \\
So we also have
\begin{equation}\label{eq:bAHd}
S_{d}\in 2\RV_{-\a,\rho}(b_{d},A_{d},H_{d}), \;\text{ with} \;\,  H_{d}(x)=c_d \,x^{-\a} \frac{(x^{\rho}-1)}{\rho}, 
\end{equation}
and $b_{d}(t)\underset{t\to\infty}{\sim} b^{*}(t)$, $A_{d} (t)\underset{t\to\infty}{\sim} A^{*}(t)$ such that $A_{d}(b_{d}(t))\underset{t\to\infty}{\sim}A(b(t))$. 
\item[]
\item[(2)]  The diversification index of $\bX$ satisfies \eqref{lim:D1minxtk}, {\it i.e.}, for $x>0$,
\begin{equation}\label{eq:DB2MRV}
 \lim_{t \to \infty} \frac{D_{1- x/t}(\bX) - K_{d} }{A_{d}(b_{d}(t))}\, =
\frac{(c_{d}-c_{1})K_{d}}{\alpha\rho} (x^{{-\rho/\a}}-1).
\end{equation}
\end{enumerate}
\end{Theorem}

\begin{proof}[Proof of Theorem \ref{thm:2mrvtoSrv}]~
\begin{enumerate}
\item[(1)]  Let  $\chi_{+}$ and $\chi_{-}$  be the positive and negative parts of the signed measure $\chi$ in its Jordan decomposition \cite{neveu:1965}.
Similarly, let $\mu_{t}^{+}$ and $\mu_{t}^{-}$  are the positive and negative parts of the signed measure $\mu_{t}$ defined in \eqref{signmeas} in its Jordan decomposition.  Definition  \eqref{signmeas} holds valid under both Assumptions \ref{assumption1} and \ref{assumption2}. Thus we have
 $\mu_{t}=\mu_{t}^{+}-\mu_{t}^{-}$  and $\chi=\chi_{+}-\chi_{-}$ where  $\mu_{t}^{+}, \mu_{t}^{-}, \chi^{+}, \chi^{-}$ are all  positive Radon measures.
Now, since $\bX\in 2\MRV_{-\alpha,\rho}(b,A,\nu,H)$, applying either \cite[Proposition 5]{resnick:2002} if Assumption \ref{assumption1} holds, or \cite[Theorem 2 ]{resnick:2002} if  Assumption \ref{assumption2}  holds, provides, as $t\to\infty$, 
$$\mu_{t}^{\pm}\conv \chi^{{\pm}}, \quad  \text{in} \quad  \M_{+}(\E). $$
Hence
\begin{align}\label{eq:2mrv_set}
\mu_{t}(\Lambda_{d}) = \frac{t \, \P\left(\frac{\bX}{b(t)} \in \Lambda_d \right) - \nu(\Lambda_d)}{A(b(t))} \;\underset{t\to\infty} {\to} \; \chi(\Lambda_d) 
\end{align}
for any relatively compact $\Lambda_d \subset \E$.  Define $b^{*}(t)= (\nu(\Gamma_d))^{1/\alpha} b(t)$. 
Then, for $x>0$,
\begin{align}
 \lim_{t\to\infty} t\, \P\left({S_d}{>xb^{*}(t)}\right)  & = \lim_{t\to\infty}  t\, \P\left(\frac{\bX}{b(t)}\in x(\nu(\Gamma_d))^{1/\alpha} \Gamma_d \right)
 						       =  \nu\left(x(\nu(\Gamma_d))^{1/\alpha} \Gamma_d\right ) = x^{-\alpha} \label{eq:sdbt}.
\end{align}
 Now, let, for $t>0$, 
$A^*(t) =A((\nu(\Gamma_d))^{-1/\alpha}t) \underset{t\to\infty} {\sim}  (\nu(\Gamma_d))^{-\rho/\alpha}\,A(t)$ (since $|A|$ is $\RV_{\rho}$). 

Let $F_{d}(x)=\P(S_{d}\le x)$. Now, for $x>0$ and $t'>0$,
\begin{align*}
\frac{\frac{\oF_{d}(t'x)}{\oF_{d}(t')}-x^{-\alpha}}{A^{*}(t{'})} & =  \frac{\frac{\P(S_{d}>b^{*}(t)x)}{\P(S_{d}>b^{*}(t))}-x^{-\alpha}}{A^{*}(b^{*}(t))} \quad \quad  \text{(replacing $t'$ by $b^{*}(t)$)} \\
  & = \frac{t\;\P(S_{d}>b^{*}(t)x)- x^{-\alpha}}{t\;\P(S_{d}>b^*(t)) A(b(t))} - x^{-\alpha} \frac{t\;\P(S_{d}>b^{*}(t))- 1}{t\;\P(S_{d}>b^*(t))A(b(t))}\\
 & = \frac{1}{t\;\P(S_{d}>b^*(t))} \times \frac{t\, \P\left(\frac{\bX}{b(t)}\in x(\nu(\Gamma_d))^{1/\alpha} \Gamma_d \right)- \nu(x(\nu(\Gamma_d))^{1/\alpha} \Gamma_d)}{A(b(t))}\\
   & \quad\quad\quad - x^{-\alpha}\frac{1}{t\;\P(S_{d}>b^*(t))}  \frac{t\, \P\left(\frac{\bX}{b(t)}\in (\nu(\Gamma_d))^{1/\alpha} \Gamma_d \right)- \nu((\nu(\Gamma_d))^{1/\alpha} \Gamma_d)}{A(b(t))}\\
   & =  \frac{1}{t\;\P(S_{d}>b^*(t))} \mu_{t}\left( x(\nu(\Gamma_d))^{1/\alpha} \Gamma_d \right) -  x^{-\alpha}\frac{1}{t\;\P(S_{d}>b^*(t))} \mu_{t}\left( (\nu(\Gamma_d))^{1/\alpha} \Gamma_d \right).
\end{align*}
Taking limit as $t'\to \infty$, and applying \eqref{eq:2mrv_set} and \eqref{eq:sdbt} in this last equation, we obtain
 \begin{align}\label{eq:chiminchi}
\lim_{t'\to \infty} \frac{\frac{\oF_{d}(t'x)}{\oF_{d}(t')}-x^{-\alpha}}{A^{*}(t{'})}  & =  \chi\left( x(\nu(\Gamma_d))^{1/\alpha} \Gamma_d \right) -  x^{-\alpha}\chi\left( (\nu(\Gamma_d))^{1/\alpha} \Gamma_d \right).
\end{align}
Since,  by assumption, $\chi(z\Gamma_{d})$ is neither a multiple of $z^{-\alpha}$, nor 0, then the same statement holds for the RHS of \eqref{eq:chiminchi}. Hence, by Remarks \ref{rk:equivResDeHaan} and \ref{rem:Hnondeg}, (or (\cite[Theorem 1]{dehaan:stadtmueller:1996})), we have $S_{d}\in 2\RV_{-\alpha,\rho}$. We need to show that  \eqref{eq:bAHd} holds with $H_d(x)=c_dx^{-\alpha}(x^{\rho}-1)/\rho$. Since $S_{d}\in 2\RV_{-\alpha,\rho}$, the RHS of \eqref{eq:chiminchi} is of the form given in \eqref{eq:def2RValt}. So we can write
\begin{align*}
\chi\left( x(\nu(\Gamma_d))^{1/\alpha} \Gamma_d \right) -  x^{-\alpha}\chi\left( (\nu(\Gamma_d))^{1/\alpha} \Gamma_d \right) = c_{a} x^{-\alpha} \frac{x^{\rho}-1}{\rho}
\end{align*}
for some $c_{a}\neq 0$,
which implies
\begin{align}\label{eq:chic1c2}
\chi\left( x(\nu(\Gamma_d))^{1/\alpha} \Gamma_d \right) = c_{a} x^{-\alpha} \frac{x^{\rho}-1}{\rho} + c_{b} x^{-\alpha}
\end{align}
where $c_{b}=\chi\left( (\nu(\Gamma_d))^{1/\alpha} \Gamma_d \right).$  We can also show (similar to \eqref{eq:chiminchi}) that
 \begin{equation*}
\lim_{t\to \infty}  \frac{t\;\P(S_{d}>b^{*}(t)x)-x^{-\alpha}}{A^{*}(b^{*}(t))} =  \chi\left( x(\nu(\Gamma_d))^{1/\alpha} \Gamma_d \right) 
     = c_{a} x^{-\alpha} \frac{x^{\rho}-1}{\rho} + c_{b} x^{-\alpha} 
\end{equation*}
by using \eqref{eq:chic1c2} in the last equality.

Therefore, by Remark \ref{rem:other2rv}, we have $S_{d}\in2\RV_{-\alpha,\rho}(b^{*},A^{*},H_{c_a,c_b})$.
If $c_{b}=0$, then 
\begin{equation}\label{eq:bAHd0}
S_{d}\in2\RV_{-\alpha,\rho}(b^{*},A^{*},H^{*}) \quad \text{where} \quad 
H^{*}(x)=  \chi\left( x(\nu(\Gamma_d))^{1/\alpha} \Gamma_d \right) = c_{a} x^{-\alpha} \frac{x^{\rho}-1}{\rho}
\end{equation}
 as per Definition \ref{def:2rv} and we define $c_d:=c_a\neq 0$. 
On the other hand, if $c_{b}\neq 0$, then, via Proposition \ref{prop:alt2rv}, we have
\begin{equation}\label{eq:bAHd1}
S_{d}\in 2\RV_{-\a,\rho}(b_{d},A_{d},H_{d}), \quad \text{ with} \;H_{d}\;\text{of the form in} \;\eqref{eq:def2RV},  
\end{equation}
$b_{d}(t)\underset{t\to\infty}{\sim} b^{*}(t)$, $A_{d} (t)\underset{t\to\infty}{\sim} A^{*}(t)$ such that $A_{d}(b_{d}(t))\underset{t\to\infty}{\sim}A(b(t))$. Here $c_d$ is obtained from $H_d(x)=c_dx^{-\alpha}(x^{\rho}-1)/\rho$ using Proposition  \ref{prop:alt2rv}. Combining \eqref{eq:bAHd0} and \eqref{eq:bAHd1} gives \eqref{eq:bAHd}. 
\item[]
\item[(2)] Since the marginal distributions are identical and $X_{1}\in 2\RV_{-\alpha,\rho}(b,A,H_1)$ as defined in \eqref{eq:def2RVinprop-bis},
we can deduce the result by combining part (1) of Theorem~\ref{thm:2mrvtoSrv}  and  part (2.a.i.)  in Theorem~\ref{thm:riskconc} (Eq. \eqref{lim:D1minxtk}) with $\displaystyle\kappa=\lim_{t\to\infty}A(b(t))/A_d(b_{d}(t))=1$. 
\vspace{-4ex}
\end{enumerate}
\end{proof}

\begin{Remark}
A sufficient condition for $S_{d}\in 2\RV$ is that $\chi(z\Gamma_{d})$ is not a multiple of $z^{-\alpha}$ 
(nor 0). It is illustrated in Example~\ref{ex:X2mrvbutSdNOT2rv}, where the condition is violated and, even though $\bX\in 2\MRV$, we see that $S_{d}\notin 2\RV$. 
\end{Remark}

\section{Examples}\label{sec:ex}

In this section, we exhibit our results using a few examples. Additionally, we provide an example where $\bX\in2\MRV$ would not imply $S_d\in2\RV$, and cases where identical univariate marginals being $2\RV$ may not imply that the random vector formed with those marginals is $2\MRV$.

\subsection{Diversification under $2\MRV$}
We develop two examples possessing $2\MRV$, one when there exists a density (Assumption \ref{assumption1}), and the other one when not (Assumption \ref{assumption2}). In both cases we can apply Theorem~\ref{thm:riskconc} to compute the asymptotic limit for the diversification index $D_{\beta}$. The dimension is restricted to $d=2$ for convenience. 

\begin{Example}{\it Pareto-Lomax marginal distribution with survival Clayton copula}
 \label{ex:parclay-lom}
~\vspace{.7ex}
\normalfont

Let us revisit Example~\ref{subsec:Exappend}.
Suppose $\bX=(X_{1}, X_{2}) \sim F$ with identical $(\alpha, 1)$-Pareto-Lomax marginal distributions, with $\alpha>1$, such that 
$$\displaystyle \overline F_1(x)=\overline F_{2}(x)=\left(1+x\right)^{-\alpha}, \quad \forall x>0.$$
Assume the dependence structure of $\bX$ to be given by a survival Clayton copula on $[0,1]^{2}$, with parameter $\theta>0$ : 
\begin{align*}
\! \P(X_1>x_1,X_2>x_2) = \left[(\overline{F}_1(x_1))^{\theta}+(\overline{F}_2(x_2))^{\theta}-1\right]^{-1/\theta} = \left[(1+x_1)^{\a \theta}+(1+x_2)^{\a \theta}-1\right]^{-1/\theta}\!\!\!.
\end{align*}
$\bullet$ 
First we verify  that $\bX \in 2\MRV$ and identify all the parameters. With $b(t)= t^{1/\alpha}-1$, we have
\begin{equation}\label{valueNu1}
t\;\P\left(\frac{\bX}{b(t)}  \in ([0,x_1]\times[0,x_2])^c\right) 
\underset{t\to\infty}{\to}  x_1^{-\a} + x_2^{-\a} - \left(x_1^{\alpha\theta} +x_2^{\alpha\theta} \right)^{-1/\theta} =:  \nu\left( ([0,x_1]\times[0,x_2])^c\right).
\end{equation}
Choosing $\displaystyle A(t)= -(t+1)^{-\min(\a\theta,1)}$, we obtain 
$$
 \lim_{t\to\infty} \frac{t\,\P\left(\frac{\bX}{b(t)} \in (([0,x_1]\times[0,x_2])^c\right) -  \nu(([0,x_1]\times[0,x_2])^c)}{A(b(t))}  = H(x_1,x_2) , \quad \text{with}
$$
{\small
\begin{equation}\label{eq:H-ex}
H(x_1,x_2) :=
\left\{ \begin{array}{lcl}
 \frac1{\theta} \left(x_1^{\alpha\theta} +x_2^{\alpha\theta} \right)^{-1-\frac1{\theta}}&\text{if}& \a\theta<1 \\
&&\\
 \alpha\left[ \left(x_1+ x_2\right)^{-(\a +1)} (x_1+x_2-1) - x_1^{-(\a +1)} (x_1-1) - x_2^{-(\a +1)}(x_2-1)   \right] & \text{if} & \a\theta=1 \\
&&\\
\alpha  \left[  \left(x_1^{\alpha\theta} +x_2^{\alpha\theta} \right)^{-1-\frac1{\theta}}\left[x_1^{\a\theta -1} (x_1-1) + x_2^{\a\theta-1} (x_2-1) \right] \right.\\
\left. \qquad \qquad\qquad \qquad\qquad-x_1^{-(\a +1)} (x_1-1) - x_2^{-(\a +1)} (x_2-1) \right]& \text{if} & \a\theta>1
\end{array}\right.
\end{equation}
}
from which we deduce that 
\begin{equation}\label{expleP-C:M2RV} 
\bX \in 2\MRV_{-\a,-1}(b,A,\nu,H) \; \text{with} \; 
\left\{\begin{array}{l}
b(t)=t^{1/\a}-1\\
A(t)= -(t+1)^{-\min(\a\theta,1)} \\
\nu  \;\text{defined in} \; \eqref{valueNu1} \\
 H \;\text{defined in} \; \eqref{eq:H-ex}.
\end{array}
\right.
\end{equation}
$\bullet$ For the next step, we  compute the density function $f$ of the distribution function $F$, as well as the density function $\lambda$ of the limit measure $\nu$, and obtain:
\begin{equation}\label{eq:density-f} 
f(x_1,x_2)=\a^2(1+\theta) \,(1+x_1)^{\a\theta-1} \,(1+x_2)^{\a\theta-1} \left((1+x_1)^{\alpha\theta} +(1+x_2)^{\alpha\theta}-1\right)^{-\frac1{\theta}-2} 
\end{equation}
and
\begin{equation}\label{eq:density-lambda} 
 \lambda(x_1,x_2) = \a^2(1+\theta)\, x_1^{\a\theta-1} \, x_2^{\a\theta-1} \left(x_1^{\alpha\theta} +x_2^{\alpha\theta} \right)^{-\frac1{\theta}-2} .
\end{equation}
%
Since $F$ has a density, we turn to Assumption \ref{assumption1}, considering for instance the case $\a\theta\ge 1$ so that $\displaystyle A(t)= -(t+1)^{-1}$. 
Checking Assumption \ref{assumption1} boils down to verifying conditions \eqref{cond4} and \eqref{cond5}. 

The following is an analysis when we have $\alpha\theta=1$; the alternative case $\alpha\theta >1$ is analogous and skipped for this part. 
Relations \eqref{eq:density-f} and  \eqref{eq:density-lambda} simplify to
\begin{equation*}
 f(x_1,x_2) = \a(\a+1) (1+x_1+x_2)^{-(\a+2)} \quad \text{and} \quad \lambda(x_1,x_2) = \a(\a+1) (x_1+ x_2)^{-(\a+2)}. 
\end{equation*}
Hence, we have, for any $\bx\in \E$, \,
$$\displaystyle 
\frac{f(t\bx)}{t^{-2}\overline{F}_1(t)} - \lambda(\bx)=
\lambda(\bx)\, t^{-1} \left (\a - \frac{2+\a}{x_1+x_2}\right)   \;\underset{t\to\infty}{\longrightarrow} \;  0.$$
Therefore, \eqref{cond4} holds and from the form of  $\frac{f(t\bx)}{t^{-2}\overline{F}_1(t)} - \lambda(\bx)$, it is clearly bounded if $\lambda(\bx)$ is, which is true for $\bx\in \aleph$. Thus uniform convergence also holds. Conditions \eqref{cond5} can be checked in a similar manner and is omitted here.\\[0.7ex]
%
$\bullet$ Since $\bX\in 2\MRV_{-\a,-1}(b,A,\nu,H)$ with $b$, $A$, $\nu$, $H$ defined in \eqref{expleP-C:M2RV}, and Assumption \ref{assumption1} holds (when $\a\theta\ge 1$), we can apply Theorem~\ref{thm:2mrvtoSrv}.  We obtain, on one hand, 
$$
S_2=X_1+X_2 \;\in\; 2\RV_{-\a,-1}(b^*,A^*,H^*), 
$$ 
where $b^*$, $A^*$ and $H^*$ are defined in \eqref{eq:bAHstar} for $d=2$, with $b(t)=t^{1/\a}-1$, $A(t)=-(t+1)^{-1}$ and $\nu$ as in \eqref{valueNu1} (computed in \eqref{eq:nugamma2}). So $b^*(t) \underset{t\to\infty}{\sim} (1+\a)^{1/\a}t^{1/\a}$, $A^{*}(b^{*}(t))=-t^{-1/\alpha}$, and
$\displaystyle H^*(x)= \chi(x(\nu(\Gamma_2))^{1/\alpha} \Gamma_2)$ ($x>0$), given that $H^*$ is non-zero and not a multiple of $x^{-\alpha}$.
 
Now, as noted in Theorem~\ref{thm:2mrvtoSrv}, \eqref{eq:bAHd}, 
we may also say that there exist $b_{2}(t)\underset{t\to\infty}{\sim} b^{*}(t)$ and $A_{2} (t)\underset{t\to\infty}{\sim} A^{*}(t)$ such that $A_{2}(b_{2}(t))\underset{t\to\infty}{\sim} A(b(t))=A^{*}(b^{*}(t))=-t^{-1/\alpha}$ and, for any $x>0$, 
$$
\frac{t\,\P(S_2/b_2(t)>x) - x^{-\alpha}}{A_2(b_2(t))} \;\underset{t\to\infty}{\to}\; H_2(x)= -c_2 \,x^{-\alpha}(x^{-1}-1).
$$

On the other hand, via Theorem~\ref{thm:2mrvtoSrv}, \eqref{eq:DB2MRV}, 
the diversification index of $\bX$ can be expressed as 
$$
 \lim_{t \to \infty} \frac{D_{1- x/t}(\bX) - K_2 }{A_2 (b_2(t))}\, =\,
\frac{(c_1-c_2)K_{2}}{\alpha} (x^{{1/\a}}-1), \quad \text{with}\; K_{2} = \frac 12\left(\frac{\nu(\Gamma_2)}{\nu(L_1)}\right)^{1/\alpha}.
$$
$\bullet$ We compute the various constants, mainly $c_1$ and $c_2$, to fully identify the results. It can be done with varying degrees of effort depending on the values of $\alpha$ and $\theta$. We restrict to the case where $\alpha\theta=1$, since the other cases would require some numerical integration, which is doable for various cases, yet not that illustrative for our purposes. \\
$\triangleright$  From Example \ref{ex:parlomax}, we have $X_{1}\in 2\RV_{-\alpha,-1} (b, A_1, H_1)$ with $A_1=-A$, $\displaystyle H_1(x) = -\a x^{-\a}(x^{-1}-1)$, hence $c_1=\a$. \\[0.5ex]
$\triangleright$ Next compute $H^*$ using $H$. 
Note that, if $\chi_{+}$ and $\chi_{+}$ are measures such that $\chi=\chi_{+}-\chi_{-}$, then
\begin{align*}
H(x_{1},x_{2}) & = \chi([(0,0),(x_{1},x_{2})]^{c})\\
 & =  \chi_{+}([(0,0),(x_{1},x_{2})]^{c}) -  \chi_{-}([(0,0),(x_{1},x_{2})]^{c})\\
 & =: \int\limits_{[(0,0),(x_{1},x_{2})]^{c}} \!\!\!\!\!\! h_{+}(x_{1},x_{2}) \;\mathrm dx_{1} \mathrm dx_{2}- 
\int\limits_{[(0,0),(x_{1},x_{2})]^{c}} \!\!\!\!\!\! h_{-}(x_{1},x_{2})\;\mathrm dx_{1} \mathrm dx_{2}.
\end{align*}
Denoting $h(x_{1},x_{2}) = h_{+}(x_{1},x_{2})-h_{-}(x_{1},x_{2})$ and differentiating $H$ defined in \eqref{eq:H-ex} for $\alpha\theta=1$ w.r.t. to $x_{1},x_{2}$ and taking care of the sign change, we obtain
\begin{equation*}
h(x_1,x_2) :=
\a(\a+1)(\a+2)( x_1+x_2)^{-(\a +3)}- \a^{2}(\a+1)( x_1+x_2)^{-(\a +2)}.
\end{equation*}
Denoting  $\displaystyle k=x\left(\nu(\Gamma_2)\right)^{1/\a}= x(1+\a)^{1/\a}$ (using \eqref{eq:nugamma2}), we have
\begin{align*}
H^*(x) &=\chi\left(k\Gamma_{2}\right)  = \int_{x_{1}+x_{2}>k} h(x_1,x_2)\, \mathrm dx_1 \mathrm dx_2\,\\
&= \a(\a+1)(\a+2) \int_{k\Gamma_{2}} ( x_1+x_2)^{-(\a +3)}\,  \mathrm dx_1 \mathrm dx_2\;- \a^2(\a+1)\int_{k\Gamma_{2}} ( x_1+x_2)^{-(\a +2)}\,  \mathrm dx_1 \mathrm dx_2 \\
 &= \,\a(\a+1)(\a+2)\int_{2k}^\infty v^{-(\a+3)}\int_0 ^v \; \mathrm du \mathrm dv \; -  \a^2(\a+1)\int_{2k}^\infty  v^{-(\a+2)} \int_0 ^v\; \mathrm du\,\mathrm dv \\
&=  \a(\a+2) \,k^{-(\a+1)} \; - \a(\a+1)\,k^{-\a} = \a x^{-\a}\left[(\a+2)(\a+1)^{-(1+1/\a)}x^{-1}-1\right],
\end{align*}
replacing $ k= x(1+\a)^{1/\a}$ in the last step.

Note that $H^*$ is of the form of \eqref{eq:def2RVHcc} with $\tilde c=-\a$ and $M=(\a+2)(\a+1)^{-(1+1/\a)}$, from which we deduce $c=M\tilde c=-\a(\a+2)(\a+1)^{-(1+1/\a)}$ and $c^*=\a\left((\a+2)(\a+1)^{-(1+1/\a)} -1\right)$ in \eqref{eq:def2RVinprop}. But we can also have the form \eqref{eq:def2RVinprop-bis} as follows:
$$
\frac{t\,\P(S_2/b_2(t)>x) - x^{-\alpha}}{-t^{-1/\a}} \;\underset{t\to\infty}{\to}\;  \a(\a+2)(\a+1)^{-(1+1/\a)}x^{-\a}(x^{-1}-1)\,,
$$
which is the same result as in Example \ref{subsec:Exappend} (see \eqref{eq:ex3.7}). Hence we can apply Theorem \ref{thm:riskconc} (as done in Example  \ref{subsec:Exappend})  to obtain the diversification benefit given by \eqref{eq:DVex3.7}.
\end{Example}
\vspace{1ex}

\begin{Example}{\it Example exhibiting asymptotic independence}.

\normalfont
Suppose $X_{1}, X_{2}$ are iid random variables with distribution function $F$ such that
\begin{align*}
\overline F(x) & =  \frac 1{2} x^{-2}(1+x^{{-1}}), \quad x\ge 1.
\end{align*}
By choosing $\displaystyle{b(t)  ={\sqrt  \frac t2} \left(1+\frac 1{\sqrt{2t}}\right)}$ and $A(t) = t^{-1}$, we observe that
\begin{align*}
t\,\P(X_{1}>xb(t)) & =  t \frac 12 x^{{-2}} \frac{2}{t}\left(1+\frac1{\sqrt {2t}}\right)^{-2} \left(1+x^{-1}{\sqrt \frac 2t} \left(1+ \frac{1}{\sqrt{2t}}\right)^{-1}\right)\\
                            & = x^{-2} + \sqrt{2}x^{-2}(x^{-1}-1) \frac 1{\sqrt t} + o\left(\frac 1{\sqrt t}\right),
\end{align*}
 and hence, \,$\displaystyle \frac{t\,\overline{F}(xb(t))-x^{-2}}{A(b(t))}  \underset{t\to\infty}{\to}  x^{-2}(x^{-1}-1)$.
 Therefore $X_{1} \in 2\RV_{2,-1}(b,A,H_{1})$ where \\ $H_{1}(x) = -c_{1}x^{-2}(x^{-1}-1)$ with $c_{1}=-1$. Moreover,  as $t\to\infty$, with $\nu_{2}(x,\infty) = x^{-\a}$, and $x>0$,
\begin{align*}
\mu_{t1}^{\pm} (x,\infty) &:= \left(\frac{t\,\P\left[\frac{X_{1}}{b(t)}\in (x,\infty) \right] -\nu_{2}(x,\infty)}{A(b(t))}\right)^{\pm}  \to  \left(x^{-2}(x^{-1}-1) \right)^{\pm}=:\chi_{1}^{\pm}(x,\infty) ,
\end{align*}
so \eqref{cond7} is satisfied.
 In this case, we know that the sum $S_2\in 2\RV$ using \cite[Theorems 3.4 and 3.5]{mao:hu:2013}; nevertheless we apply Theorem \ref{thm:2mrvtoSrv} for illustration. 

To check that $\bX=(X_{1},X_{2})\in \MRV$, take a set of the form $[0,(x_{1},x_{2})]^{c}$ for $x_{1}>0,x_{2}>0$, and observe that 
\begin{align*}
t\,\P\left(\frac{\bX}{b(t)} \in [0,(x_{1},x_{2})]^{c}\right)  \underset{t\to\infty}{\to} \frac{1}{x_{1}^{2}} +\frac{1}{x_{2}^{2}}=: \nu([0,(x_{1},x_{2})]^{c}).
\end{align*}
Hence Assumption \ref{assumption2} is satisfied and we may apply Theorem \ref{thm:2mrvtoSrv} once we have checked that \linebreak $\bX \in 2\MRV$. We can write
\begin{align*}
& \frac{t\,\P\left(\frac{\bX}{b(t)} \in [0,(x_{1},x_{2})]^{c}\right) -\nu([0,(x_{1},x_{2})]^{c}) }{A(b(t))} \\ 
                                 & = \mu_{t1}(x_{1},\infty)^{+} -  \mu_{t1}(x_{1},\infty)^{-} +  \mu_{t1}(x_{2},\infty)^{+}-  \mu_{t1}(x_{2},\infty)^{-} + o(1)\\
                                 &  \underset{t\to\infty}{\to} \chi_1^{+}(x_{1},\infty) - \chi_1^{-}(x_{1},\infty) + \chi_1^{+}(x_{2},\infty)- \chi_1^{-}(x_{2},\infty)= x_{1}^{-2}(x_{1}^{-1}-1)  + x_{2}^{-2}(x_{2}^{-1}-1).
\end{align*}
Hence $\bX \in 2\MRV_{-2,-1}(b,A, \nu, H)$ with $H(x_{1},x_{2}):=x_{1}^{-2}(x_{1}^{-1}-1)  + x_{2}^{-2}(x_{2}^{-1}-1)$. 
Note that $\displaystyle \nu(\Gamma_{2}) = \nu(\{(x_{1},x_{2}): x_{1}+x_{2}>1\}) = 1+1 =2$.
 Defining $\chi([0,(x_{1},x_{2})]^{c}) = H(x_{1},x_{2})$ for $x_1>0$, $x_{2}>0$, we see that $\chi$ also concentrates on the axes. Hence, for $k>0$,
 \begin{align*}
 \chi(k\Gamma_{2}) = H(k,k) = 2k^{-2}(k^{-1}-1). \end{align*}
Therefore from Theorem \ref{thm:2mrvtoSrv}(1), with $b^{*}(t) = \sqrt{2}\,b(t) , A^{*}(t) = \sqrt{2}\,A(t)$ and $H^{*}(x)=\chi(\sqrt{2}\,x\Gamma_{2})$, we have $S_{2}\in 2\RV_{-2,-1}(b^{*},A^{*},H^{*})$ where
\begin{align*}
H^{*}(x)=\chi(\sqrt{2}\,x\Gamma_{2}) =x^{-2}\left(\frac{1}{\sqrt{2}}x^{-1}-1\right) = \frac{c}{\rho}x^{-\a}(M x^{\rho}-1),\;\text{with}\; c=-1,\,\rho=-1,\,\a=2,\,M=\frac1{\sqrt{2}}.
\end{align*}
Using Remark~\ref{rem:other2rv} and Proposition~\ref{prop:alt2rv}, we have $S_{2}\in 2\RV_{2,-1}(b_{2},A_{2},H_{2})$ where
\begin{align*}
b_{2}(t) &= (1+A^{*}(b^{*}(t)))b^{*}(t)=\sqrt{t}+\sqrt{2}+\frac{1}{\sqrt{2}}\,\\
A_{2}(t) &= A^{*}(b^{*}(b_{2}^{\la}(t))) = \frac{\sqrt{2}}{t-\sqrt{2}-1/\sqrt2}=\left(\frac{t}{\sqrt2} - \frac32\right)^{-1},
\end{align*}
and 
$\displaystyle H_{2}(x) =  \frac{c_{2}}{\rho}x^{-\a}(x^{\rho}-1)$, with $c_{2}=c\,M= -1/\sqrt{2}$.

Hence using Theorem \ref{thm:2mrvtoSrv}, \eqref{eq:DB2MRV}, where $K_{2}=\frac 12\left(\frac{\nu(\Gamma_2)}{\nu(L_1)}\right)^{1/\alpha} \!\!\!=1/\sqrt 2$, we obtain, for $x>0$,
$$
 \lim_{t \to \infty} \frac{D_{1- x/t}(\bX) - K_2 }{A_2 (b_2(t))}\, =\,
\frac{(c_2-c_1)K_{2}}{\alpha\rho} (x^{{-\rho/\a}}-1) = -\frac{\sqrt{2}-1}{4} (\sqrt{x}-1).
$$
\end{Example}

\subsection{Does $\bX\in 2\MRV$ necessarily imply that $S_{d}\in2\RV$?}

We provide an example where second order multivariate regular variation exists, yet the sum does not possess univariate second order regular variation. This justifies the condition imposed on the limit measure $\chi$ in Theorem \ref{thm:2mrvtoSrv}, namely $\displaystyle \chi(z\Gamma_d)$ is neither multiple of $z^{-\alpha}$, nor 0.
\begin{Example}\label{ex:X2mrvbutSdNOT2rv}
\normalfont

 Consider random variables $W, Z$ where 
 \begin{align*}
 \P(W>w) & = \frac 12 w^{-1}+ \frac 12 w^{-2}, \quad w>1,\\
 \P(Z>z) & =  2 z^{-1} - z^{-2}, \quad z>1.
 \end{align*}
Let $B_{1}, B_{2}$ be Bernoulli random variables with  $\P(B_{1}=1)= 8/9 = 1-\P(B_{1}=0)$ and \linebreak $\P(B_{2}=1)= 1/2 = \P(B_{2}=0)$. Also assume that $W,Z, B_{1}, B_{2}$ are mutually independent. \\ Now define the vector: 
\[\bX= (X_{1}, X_{2})\quad \text{with} \quad X_1:= B_{1}B_2 W + (1-B_{1})Z \quad  \text{and}  \quad 
X_2:= B_{1} (1-B_2) W + (1-B_{1})Z.\]

For any $\bx=(x_1,x_2)$ with $x_{1}>1, x_{2}> 1$,
\begin{align*}
\P(\bX \in [\bzero,(x_{1},x_{2})]^{c}) & = \frac 49\cdot \P(W>x_{1}) + \frac 49\cdot \P(W>x_{2}) + \frac 19 \cdot\P(Z>x_{1} \wedge x_{2})\\
                                 & = \frac 29 \left(x_{1}^{-1}+x_{2}^{-1}+(x_{1}\wedge x_{2})^{-1}\right) +  \frac 19 \left(2x_{1}^{-2}+2x_{2}^{-2}-(x_{1}\wedge x_{2})^{-2}\right)
\end{align*}
where $x_{1}\wedge x_{2} = \min(x_{1},x_{2})$. Therefore, for $\bx\in (0,\infty)\times(0,\infty)$,
\begin{align*}
\lim_{t\to\infty} t\P\left(\frac{\bX}t \in [\bzero,\bx]^{c}\right) & = \frac 29 \left(x_{1}^{-1}+x_{2}^{-1}+(x_{1}\wedge x_{2})^{-1}\right) =:\nu([\bzero,\bx]^{c}).
\end{align*}
Moreover, 
\begin{align}\label{eq:H-X2mrvbutSdNOT2rv}
 \lim_{t\to\infty} \frac{t\P\left(\frac{\bX}t \in [\bzero,\bx]^{c}\right) -\nu([\bzero,\bx]^{c})}{t^{-1}} &=    \frac 19 \left(2x_{1}^{-2}+2x_{2}^{-2}-(x_{1}\wedge x_{2})^{-2}\right) \nonumber\\
&=\frac 19 \left((x_{1}\wedge x_{2})^{-2}+2\,(x_{1} \vee x_{2})^{-2}\right)
=: H(\bx),
\end{align}
where $x_{1}\vee x_{2} = \max(x_{1},x_{2})$. 
Hence \;$\bX\in 2\MRV_{-1,-1}(b,A,\nu,H)$ with $b(t)=t$ and $A(t)=t^{-1}$. 

Now, with $S_2=X_1+X_2$, we can write, for $x>2$,
$$
\P(S_2>x) = \P(B_1 W+2(1-B_1)Z >x)
=\frac 89\cdot \P(W>x) + \frac 19 \cdot\P\left(Z>\frac x2\right)
=  \frac 89 \,x^{-1}.
$$                                
Therefore $S_2\in \RV_{-1}$ but $S_2$ clearly does not possess second order regular variation. 

Defining $\chi([\bzero,\bx]^c):=H(\bx)$, with $H$ defined in \eqref{eq:H-X2mrvbutSdNOT2rv}, we can check that in such a case, $\chi(z\Gamma_2)=0$ for all $z>0$. Hence  Theorem \ref{thm:2mrvtoSrv} cannot be applied as one of the assumptions does not hold.

\end{Example}

\subsection{Parameter stability in second order regular variation}

In Section \ref{subsec:aggr} we discuss conditions under which assuming a vector $\bX=(X_{1},\ldots,X_{d})\in2\MRV_{{-\alpha,\rho}}$ with identical marginals $X_1\in 2\RV_{-\alpha,\rho}$ would imply that $S_{d}=X_{1}+\ldots+ X_{d}\in 2\RV_{-\alpha,\rho}$ with $\alpha>0$ and $\rho\le 0$. One may ask here, if the margins $X_{1}\in 2\RV_{{-\a,\rho}}$, with some nice dependence structure on $\bX$ (like independence), will this imply $S_d \in 2\RV_{{-\a,\rho}}$? Alternatively, we may also ask if $X_1\in 2\RV, \bX\in 2\MRV$ and $S_d\in2\RV$, does this necessitate that all of them have the same parameters of regular variation $(-\alpha,\rho)$? Assuming identical marginals, the first order parameter $-\alpha$ clearly remains the same for $X_1,\bX$ and $S_d$.

%
In the following example, we observe that actually a variety of possibilities for the second order parameter exist for $S_d\in 2\RV$ as well as for $\bX\in 2\MRV$. The example also provides a justification for the assumption $\bX\in 2\MRV_{-\alpha,\rho}$, on top of the marginal assumptions of $X_1 \in 2\RV_{-\alpha,\rho}$ (via Assumptions \ref{assumption1} or \ref{assumption2}) in Theorem \ref{thm:2mrvtoSrv}.

\begin{Example}\label{ex:1}
\normalfont
Suppose $X_{1}, X_{2}$ are iid random variables with distribution function $F$ such that
\begin{align*}
\overline F(x) & = \frac{1}{2} x^{-\alpha}(1+x^{{\rho}}), \quad x\ge 1,
\end{align*}
where $\alpha>0, \rho <0$. This family of distributions is often called the Hall-Welsh class of heavy-tailed distributions.
 Clearly,  $X_{1}\in 2\RV_{-\alpha,\rho}(b,A)$ where $b(t)=t^{{1/\a}}$ and  $A(t)=t^{{\rho}}$. Using \cite[Theorem 3.5]{mao:hu:2013},  we know that the sum $S_{d} \in 2\RV_{-\alpha,\gamma}$ where $\gamma=-\min(1,\alpha,-\rho)$. 
 Let us check if $\bX=(X_{1},X_{2})\in 2\MRV$. Take a set of the form $[0,(x_{1},x_{2})]^{c}$ for $x_{1}>0,x_{2}>0$, and observe that
\begin{align*}
t\,\P\left(\frac{\bX}{t^{1/\a}} \in [0,(x_{1},x_{2})]^{c}\right)  \underset{t\to\infty}{\to} \frac 12\left(\frac{1}{x_{1}^{\a}} +\frac{1}{x_{2}^{\a}}\right) =: \nu([0,(x_{1},x_{2})]^{c}).
\end{align*}

Define
\begin{align}\label{eq:imp}
 H(x_{1},x_{2},t)  & = \frac{t\,\P\left(\frac{\bX}{t^{1/\a}} \in [0,(x_{1},x_{2})]^{c}\right) - \frac 12\left(\frac{1}{x_{1}^{\a}} +\frac{1}{x_{2}^{\a}}\right) }{A^{*}(t^{1/\alpha})} 
\end{align}
where $A^{*}(t) = t^{\max(\rho,-\alpha)}$. Therefore
\begin{align*}
 H(x_{1},x_{2},t)  & = \left\{ \begin{array}{ll}   \frac 12 (x_{1}^{-\a+\rho}+x_{2}^{-\a+\rho}) -\frac{1}4t^{-1-\rho/\a}x_1^{-\a}x_2^{-\a}(1+t^{\rho/\a}x_{1}^{\rho})(1+t^{\rho/\a}x_{2}^{\rho}), &  \quad\rho+\alpha\ge 0\\[0.7em]
 \frac {1}2  t^{1+\rho/\a} (x_{1}^{-\a+\rho}+x_{2}^{-\a+\rho}) -\frac14x_1^{-\a}x_2^{-\a}(1+t^{\rho/\a}x_{1}^{\rho})(1+t^{\rho/\a}x_{2}^{\rho}), &  \quad\rho+\alpha< 0.
 \end{array} \right.
\end{align*}
Now, we have
\[ \lim_{t\to\infty} H(x_{1},x_{2},t) =
  \begin{cases}
    \frac12 (x_{1}^{-\a+\rho}+x_{2}^{-\a+\rho}),      & \quad \text{if }  \rho+\a>0,\\
     \frac12(x_{1}^{-2\a}+x_{2}^{-2\a})   -\frac14 x_1^{-\a}x_2^{-\a},  & \quad \text{if }  \rho+\a=0,\\
      -\frac14x_1^{-\a}x_2^{-\a} & \quad \text{if }  \rho+\a<0.
  \end{cases}
\]
Therefore, we have $\bX\in 2\MRV_{-\a,\rho}$ if $\alpha+\rho\ge 0$, and $\bX\in 2\MRV_{-\a,-\a}$ if $\alpha+\rho< 0$. 
We can check that no other choice of $A(\cdot)$ (up to equivalent tail behavior) provides a finite limit for \eqref{eq:imp} as $t\to\infty$. 

This means that we may have a variety of indices of second order regular variation appearing together. For example, if $1<\alpha<-\rho$, then $\gamma=-\min(1,\alpha,-\rho) =\rho$ and $\alpha+\rho<0$, which implies:
\[X_{1}\in 2\RV_{-\alpha,\rho}, \quad\quad S_d\in2\RV_{-\alpha,-1}, \quad\quad   \bX\in 2\MRV_{-\alpha,-\alpha}.\]
Note that, in this case, although Assumption \ref{assumption2} is satisfied, we do not have $\bX\in 2\MRV_{-\alpha,\rho}$ and hence Theorem \ref{thm:2mrvtoSrv} cannot be used.

\end{Example}

\subsection{Other constructions leading to second order regular variation for sums}

We  have seen that the  multivariate second order regular variation conditions discussed in Section \ref{subsec:aggr} provide a class of examples for distributions whose sums are also $2\RV$. But other constructions are possible too, as shown in the following example. 

\begin{Example}\label{ex:cdtional}{\it Conditional Independence.}
\normalfont
A possible way for obtaining dependent random variables whose sum will admit $2\RV$ is to  make them  conditionally independent. Suppose $\bX=(X_{1},\ldots, X_{d})\in \R^{d}_{+}$ is a random vector and there exists a latent random variable $\Theta$ such that  $(X_i | \Theta=\theta)$, $i=1,\cdots,d$, are conditional independent and identically distributed random variables with 
$$(X_{i}|\Theta=\theta) \sim 2\RV_{-\alpha,\rho}(b_{\theta}, A_{\theta},H_{\theta})$$ 
where $\alpha>0$ and $\rho \le 0$. Let $S_d=\sum_{i=1}^d X_i$. Then using  \cite[Theorem 3.5]{mao:hu:2013}, we have that the conditional sum $$(S_d|\Theta=\theta) \in 2\RV_{-\alpha,\rho}(\bar{b}_{\theta}, \bar{A}_{\theta},\bar{H}_{\theta})$$ with appropriate $\bar b_{\theta}$, $\bar A_{\theta}$ and $\bar H_{\theta}$. Now, under certain choices of $\alpha, \rho$ and mild integrability conditions, we can show that the sum $S_d\in 2\RV_{-\alpha,\rho}$.  \\

As an example, consider $\bY=(Y_{1},Y_{2}) = \Theta(X_{1},X_{2})$ where $(X_{1},X_{2})$ has the same distribution as 
in Example  \ref{ex:1} with $\rho>-\alpha$, and 
$\Theta$ is an independent random variable of $\bX$ such that $\E[|\Theta|^{-\rho}]<\infty$. Recall that we found $X_{1}\in 2\RV_{-\alpha,\rho}(b,A)$ where $b(t)=t^{{1/\a}}$ and  $A(t)=t^{{\rho}}$.   We can check that $$(Y_{1}|\Theta=\theta)\in 2\RV_{-\alpha,\rho} (b_{\theta}, A_{\theta})$$ where 
$b_{\theta}= b(t)$ and $A_{\theta}(t)=\theta^{-\rho}A(t)$. Hence, using  \cite[Theorem 3.5]{mao:hu:2013}, we have $(S_2|\Theta=\theta) \in 2\RV_{-\alpha,\rho}(b, A_{\theta})$.   Now, using the fact that  $\E(|\Theta|^{-\rho})<\infty$ and writing $\displaystyle F_{S_2}(x)=\int F_{S_2\mid\Theta=\theta}\left(x\right) dF_{\Theta}\left(\theta\right)$, we can verify  that  $S_2 \in 2\RV_{-\alpha,\rho}(b, A)$. 
\end{Example}

\section{Conclusion}\label{sec:concl}

The motivation for this work has been to study diversification benefits in a portfolio of heavy-tailed risk factors. In doing so, we can highlight two main contributions of this paper. First, we find  the convergence rate of the diversification benefit for $\VaR$ as the level tends to 1, assuming second order regular variation for the portfolio. Secondly, we explore in detail the relationship between second order regular variation of a vector, its marginal components and their sum. Although the assumptions imposed in our results are often sufficient conditions, we exhibit the importance of these assumptions via counterexamples.

A few questions still remain open. For instance, a characterization of multivariate second order regular variation in terms of linear combination of its marginals akin to a Cram\'er-Wold Theorem is yet to be discovered. It may also be interesting to study the effects of the related concept of hidden regular variation on diversification. We intend to explore these directions of research in the near future.

\section*{Acknowledgement} 
Both authors are grateful to the referees, including the associate editor, for their insightful reviews of the manuscript and many helpful suggestions.


\bibliographystyle{spmpsci} 
\bibliography{bibvarhrv}

\section*{APPENDIX}

The diversification property relating the marginal risks to the aggregate risk in Theorem~\ref{thm:riskconc}(2.a), can be easily extended to tail equivalent risks. We provide the result in the following.
\begin{Theorem}\label{theo:appendix}
Let $X\in 2\RV_{-\alpha,\rho_{ X}}(b_{X},A_{ X},H_{ X})$ and $Y\in 2\RV_{-\alpha,\rho_{ Y}}(b_{ Y},A_{ Y},H_{ Y})$ with $\alpha>0, \rho_X<0,\rho_Y<0$ and 
\[H_{ X} (x):= \frac{c_{ X}}{\rho_{ X}}x^{-\alpha}(x^{\rho_{ X}}-1), x>0, \quad\quad H_{ Y}(x):= \frac{c_{ Y }}{\rho_{ Y}}x^{-\alpha}(x^{\rho_{ Y}}-1), x>0.\]
Assume that $X$ and $Y$ are tail equivalent risks, meaning that $\displaystyle 0<\lim_{t\to\infty} b_X(t)/b_Y(t)<\infty$, and define $\displaystyle \tilde{K}_d:=\lim_{t\to\infty} \frac{b_{ Y}(t)}{d\,b_{ X}(t)}$. Then the following hold.
\begin{enumerate}
\item[i.] If $\displaystyle{\lim_{t\to\infty}A_{ X}(b_{ X}(t))/A_{ Y}(b_{ Y}(t)) =\kappa\in \R}$, then, for any $x>0$, we have
\begin{align*}
  \lim_{t \to \infty} \frac1{ A_{ Y}(b_{ Y}(t))}\times\left(\frac{\VaR_{1-x/t}(Y)}{d\,\VaR_{1-x/t}(X)} - \tilde{K}_d\right)  
\, =\,
\frac{(c_{ Y}-\kappa c_{ X})\tilde{K}_{d}}{\alpha\rho_{ Y}} (x^{{-\rho_{ Y}/\a}}-1).
\end{align*}
\item[ii.] If $\displaystyle{\lim_{t\to\infty}A_{ Y}(b_{ Y}(t))/A_{ X}(b_{ X}(t)) =0}$, then, for any $x>0$, we have
\begin{align*}
  \lim_{t \to \infty} \frac1{ A_{ X}(b_{ X}(t))}\times\left(\frac{\VaR_{1-x/t}(Y)}{d\,\VaR_{1-x/t}(X)} - \tilde{K}_d\right) 
\, =\,
-\frac{c_{ X}\tilde{K}_{d}}{\alpha\rho_{ X}} (x^{{-\rho_{ X}/\a}}-1).
\end{align*}
\item[]
\end{enumerate}
\end{Theorem}
\begin{proof}
The proof of Theorem~\ref{theo:appendix} is the same as that of Theorem~\ref{thm:riskconc}(2.a) and can be obtained by replacing $X_1$ by $X$, $S_d$ by $Y$ (with the corresponding parameters for the $2\RV$ property), $D_\beta$ by $\displaystyle \frac{\VaR_\beta(Y)}{d\,\VaR_{\beta}(X)}$, and $K_d$ by $\displaystyle \tilde{K}_d$ in the proof of Theorem~\ref{thm:riskconc}(2.a).
\end{proof}

\end{document}